\documentclass[12pt, a4paper]{article}
\usepackage[left=1in,right=1in,top=1.0in,bottom=1in]{geometry}

\usepackage{amsmath,amssymb,amsthm,amsfonts,afterpage,hyperref}
\usepackage{amscd,subeqnarray}

\usepackage{graphicx}
\usepackage{caption}
\usepackage{subcaption}
\usepackage{authblk}
\usepackage{tikz}

\usepackage{array}
\usepackage{wrapfig}
\usepackage{color}
\usepackage{ragged2e}
\usepackage{stmaryrd}
\usepackage{siunitx}
\usepackage{multirow}
\usepackage{xcolor,cancel}
\usepackage{boldfonts}
\usepackage{symbols}
\usepackage{footmisc}
\usepackage{float}
\usepackage{setspace}
\usepackage{bm}
\usepackage{booktabs}
\usepackage{tikz, xcolor}
\usepackage{soul}

\usepackage[T1]{fontenc}
\usepackage[utf8]{inputenc}

\usepackage[linesnumbered,ruled,vlined]{algorithm2e}
\SetKwInput{KwInput}{Input}                
\SetKwInput{KwOutput}{Output} 
\usepackage{algpseudocode}

\usepackage[sort, numbers]{natbib}
\setcitestyle{square}

\usetikzlibrary{shapes,arrows, positioning,calc}	
\usetikzlibrary{arrows,snakes,backgrounds}
\usepackage{mathrsfs}

\newcommand{\bfa}[1]{\boldsymbol{#1}} 			
			
\newcommand{\bfeps}{\boldsymbol{\epsilon}}

\newcommand{\Sym}{\text{Sym}}   			%
\newcommand{\tr}{\text{tr}}       				%

\DeclareMathAlphabet{\mathpzc}{OT1}{pzc}{m}{it}

\newcommand{\bfn}{\boldsymbol{n}}	
	
\newcommand{\bfu}{\boldsymbol{u}}	
	
\newcommand{\bfv}{\boldsymbol{v}}	
\newcommand{\bfw}{\boldsymbol{w}}

\newcommand{\bfB}{\boldsymbol{B}}	
\newcommand{\bfx}{\boldsymbol{x}}

\newcommand{\bfT}{\boldsymbol{T}}		
\newcommand{\bfI}{\boldsymbol{I}}	 
\newcommand{\bfzero}{\boldsymbol{0}}

\newcommand{\bff}{\boldsymbol{f}}	
\newcommand{\bfg}{\boldsymbol{g}}	

\newtheorem{theorem}{Theorem}[section]

\newtheorem{lemma}[theorem]{Lemma}
\newtheorem{remark}{Remark}

\newtheorem*{cwf}{Continuous weak formulation}
\newtheorem*{dwf}{Discrete weak formulation}

\providecommand{\keywords}[1]
{
  \small	
  \textbf{\textit{Keywords---}} #1
}

\title{A finite element model to analyze crack-tip fields in a transversely isotropic strain-limiting elastic solid}

\author[1]{Saugata Ghosh\thanks{saugata.ghosh01@utrgv.edu}}
\author[2]{Dambaru Bhatta\thanks{dambaru.bhatta@utrgv.edu}}
\author[3]{S. M. Mallikarjunaiah\thanks{m.muddamallappa@tamucc.edu}\thanks{corresponding author}}
\affil[1,2]{School of Mathematical \& Statistical Sciences,
The University of Texas - Rio Grande Valley, 
Edinburg, Texas 78539, USA}
\affil[3]{Department of Mathematics \& Statistics,
Texas A\&M University-Corpus Christi, 
Corpus Christi, Texas 78412-5825, USA}
\date{}

\begin{document}

\maketitle
	    
\begin{abstract}
 
This paper presents a finite element model for the analysis of crack-tip fields in a transversely isotropic strain-limiting elastic body. A nonlinear constitutive relationship between stress and linearized strain characterizes the material response. This algebraically nonlinear relationship is critical as it mitigates the physically inconsistent strain singularities that arise at crack tips.  These strain-limiting relationships ensure that strains remain bounded near the crack tip, representing a significant advancement in the formulation of boundary value problems (BVPs) within the context of first-order approximate constitutive models. For a transversely isotropic elastic material containing a crack, the equilibrium equation, derived from the balance of linear momentum under a specified nonlinear constitutive relation, is shown to reduce to a second-order, vector-valued, quasilinear elliptic BVP. A robust numerical method is introduced, integrating Picard-type linearization with a continuous Galerkin-type finite element procedure for spatial discretization. Numerical results, obtained for tensile loading conditions and two distinct material fiber orientations, illustrate that the evolution of crack-tip strains occurs significantly slower than that of the normalized stresses. However, the strain-energy density is most pronounced near the crack tip, consistent with observations from linearized elasticity theory.  It is demonstrated that the framework investigated herein can serve as a basis for formulating physically meaningful and mathematically well-defined BVPs, which are essential for exploring crack evolution, damage, nucleation, and failure in anisotropic strain-limiting elastic materials.  
 
 \end{abstract}
	 
\noindent \keywords{Finite element method; Plane-strain fracture; Transversely isotropic; Strain-limiting theories of elasticity}

\section{Introduction}

The analysis of crack tip fields within transversely isotropic elastic bodies presents a considerable challenge due to the inherent anisotropy of the material, which complicates the distributions of stress and displacement in comparison to isotropic materials \cite{anderson2017fracture,broberg1999}. This anisotropy, frequently observed in composites, timber, and geological formations, necessitates specialized numerical techniques to accurately capture the complex stress singularities and deformation patterns at the crack tips \cite{chiang2004some}. Cracks within transversely isotropic elastic bodies are of critical importance across various engineering disciplines since these materials, commonly found in composites, timber, and layered geological formations, exhibit mechanical properties that are dependent on direction, significantly influencing fracture behavior \cite{kachanov1993elastic}. A comprehensive understanding of the initiation and propagation of cracks within these anisotropic structures is paramount for ensuring structural integrity and preventing catastrophic failures. The intrinsic anisotropy complicates stress distribution around crack tips, resulting in unique challenges in fracture mechanics when compared to isotropic materials. The precise prediction of crack growth and failure in these materials necessitates a thorough understanding of the stress intensity factors and energy release rates, which are directly influenced by the material's anisotropy and the applied loading conditions. Therefore, investigating crack-tip fields in transversely isotropic bodies is essential for developing reliable design criteria and predictive models for various applications, ranging from aerospace and automotive components to civil infrastructure and energy exploration.

The precise quantification of stress and strain fields surrounding geometric discontinuities, such as notches, slits, holes, and damage inclusions, represents a fundamental problem of enduring significance in both engineering practice and theoretical solid mechanics. Traditionally, analyses of these stress concentrations have been predicated upon the constitutive assumptions of linearized elasticity theory, as established in foundational works \cite{Inglis1913,lin1980singular,love2013treatise,murakami1993stress}. A well-recognized deficiency of this classical approach is its inherent prediction of physically untenable, unbounded strain singularities at discontinuity tips, a direct consequence of its first-order linear approximation of finite deformation. Consequently, a substantial corpus of research has been directed towards the development of refined constitutive models aimed at achieving more physically congruent representations of material response \cite{gurtin1975,sendova2010,MalliPhD2015,ferguson2015,zemlyanova2012,WaltonMalli2016,rajagopal2011modeling,gou2015modeling}. Nevertheless, a persistent challenge lies in reconciling the need for enhanced model fidelity with computational efficiency and experimental validation \cite{broberg1999}, as many proposed model augmentations impose significant computational burdens or present substantial hurdles in empirical verification. Moreover, applying linear elastic fracture mechanics (LEFM) to crack initiation and propagation modeling is subject to intrinsic limitations that necessitate careful consideration. Beyond the established strain singularity at the crack tip, these limitations include predicting a physically implausible blunt crack-opening profile and the potential for crack-face interpenetration, particularly within bimaterial interfaces. Notably, the issue of crack-tip singularity remains unresolved even when utilizing nonlinear elasticity frameworks, as exemplified in \cite{knowles1983} and the bell constraint model in \cite{tarantino1997}. Therefore, a salient question emerges regarding the capacity of algebraic nonlinear models to effectively regulate the crack-tip strain singularity, even in the presence of singular stresses. 

A generalized framework for elasticity, extending beyond the classical Cauchy and Green formulations, has been developed by Rajagopal and collaborators in a series of publications \cite{rajagopal2003implicit,rajagopal2007elasticity,rajagopal2007response,rajagopal2009class,rajagopal2011non,rajagopal2011conspectus,rajagopal2014nonlinear,rajagopal2016novel,rajagopal2018note}. This body of work, \textit{Rajagopal's theory of elasticity}, introduces implicit constitutive models grounded in a robust thermodynamic foundation. The response of an elastic body, defined herein as a material incapable of energy dissipation, is effectively characterized by implicit relations between the Cauchy stress and deformation gradient tensors \cite{bustamante2018nonlinear}. A salient feature of Rajagopal's approach is the potential to derive a hierarchical structure of 'explicit' nonlinear relationships expressing linearized strain as a nonlinear function of Cauchy stress. Notably, a distinct subclass of these implicit models facilitates the representation of linearized strain with a uniformly bounded function throughout the material domain, even under conditions of substantial stress. This 'limiting strain' property renders these models particularly well-suited for investigating crack and fracture behavior in brittle materials \cite{rajagopal2011modeling,gou2015modeling,Mallikarjunaiah2015,MalliPhD2015}, including the potential for extension to quasi-static and dynamic crack evolution analyses. Utilizing these strain-limiting models, numerous studies have revisited classical elasticity problems \cite{kulvait2013,rajagopal2018bodies,bulivcek2014elastic,erbay2015traveling,itou2018states,zhu2016nonlinear,csengul2018viscoelasticity,itou2017contacting,yoon2022CNSNS,yoon2022MMS}. Strain-limiting constitutive models offer a versatile framework for elucidating the mechanical behavior of a broad spectrum of materials. This versatility is particularly advantageous in the analysis of crack and fracture phenomena. Recent investigations, as reported in \cite{lee2022finite,yoon2021quasi}, have demonstrated that the formulation of quasi-static crack evolution problems within the strain-limiting theoretical framework yields a diverse array of complex crack patterns, notably including the observation of increased crack-tip propagation velocities.  

This study investigates the behavior of a singular crack embedded within a transversely isotropic solid, employing a constitutive formulation derived from Rajagopal's theory of elasticity. A specialized constitutive relationship, tailored to represent the stress-strain response of transversely isotropic materials accurately, is developed. The combination of the linear momentum balance and the algebraically nonlinear constitutive equation results in a vector-valued, quasi-linear elliptic boundary value problem. A finite element-based numerical methodology is employed to approximate the solution due to the inherent intractability of analytical solutions for such nonlinear partial differential equations. The finite element method, renowned for accurately capturing crack-tip fields in elastic materials, provides a flexible framework for domain discretization and numerical solution of the governing partial differential equations. To address the inherent nonlinearities of the system, Picard's iterative algorithm is implemented, and convergence of the numerical solution is demonstrated through the progressive reduction of the residual at each iteration. Several intriguing results regarding stress concentration, the slow growth of cracked-tip strains, and the decrease of strain-energy density have been reported for a single crack subjected to tensile loading. This study is fundamental and may be further extended in several directions, including the examination of thermo-elastic static and quasi-static cracks, as well as dynamic crack propagation in transversely isotropic materials.  

The organization of this article is as follows: In Section~\ref{math_formulation}, the implicit theory is presented, and the derivation of the nonlinear constitutive relation is detailed. A mathematical model, describing a static crack in a transversely isotropic solid under tensile loading, is presented in Section~\ref{BVP}, where the existence of a unique solution to the weak formulation is also demonstrated. A numerical method, based on continuous Galerkin-type finite elements coupled with Picard's iterative algorithm, is presented in Section~\ref{fem}. A detailed analysis of the numerical solution and the effects of various parameters are presented in Section~\ref{rd}. Finally, a conclusion is presented in the concluding section of the paper. 

\section{Mathematical formulation}\label{math_formulation}
Let $\mathcal{D}$ be a bounded domain in $\mathbb{R}^2$, and it is assumed to be occupied by a transversely isotropic elastic body.  The boundary $\partial \mathcal{D} = \overline{\Gamma_N} \cup \overline{\Gamma_D}$ is assumed to Lipschitz, where $\Gamma_N$ is the Neumann boundary and $\Gamma_D \neq \emptyset$ is the Dirichlet boundary.  Let $\bfa{n}$ be the outward unit normal to $\mathcal{D}$. Let $\Gamma_c$ be a 1-dimensional manifold, completely contained in  $\mathcal{D}$, splitting the domain $\mathcal{D}$ into two parts. Let $\Sym(\mathbb{R}^{2 \times 2})$ be a vector space of symmetric $2 \times 2$ tensors with inner product $\bfa{A} \colon \bfa{B} = \sum_{i,\, j = 1}^2  \, \bfa{A}_{ij} \, \bfa{B}_{ij} $ and the associated \textit{induced norm} $ \| \bfa{A} \| = \sqrt{\bfa{A} \colon \bfa{A}}$. Let $\bfa{u} \colon \mathcal{D} \to \mathbb{R}^2$ be a displacement vector defined on spatial points $\bfa{x}$ in the deformed configuration. The points $\bfa{X}$ are assumed to be in reference configuration. Hence, the displacement vector is also defined as $\bfa{u} := \bfa{x} - \bfa{X}$. 

Let $L^{p}(\mathcal{D})$ denote the space of all \textit{Lebesgue integrable functions} for $p\in[1, \infty)$, and $\left( \cdot, \; \cdot \right)$ and $\| \cdot \|$ {respectively} denote the usual inner-product and the norm of functions defined on the usual {space of square-integrable functions} $L^2(\mathcal{D})$. Let $W^{k, \, p}(\mathcal{D})$ is the vector space of functions that are in $L^{p}(\mathcal{D})$ with derivaties of order up to $k$. Let $C^{m}(\mathcal{D}), \; m \in \mathbb{N}_0$ be the space of continuous functions on $\mathcal{D}$. Further, let $H^{1}(\mathcal{D})$ denote the classical \textit{Sobolev space} \cite{ciarlet2002finite,evans1998partial}:
\begin{subequations}
\begin{align}
H^{1}(\mathcal{D}) &:= \left\{ v \in L^{2}(\mathcal{D}) \; \colon \; Dv \in L^{2}(\mathcal{D}) \right\},  \label{1a} \\
{H^1_0(\mathcal{D})} & := \left\{ v \in H^1(\mathcal{D}) \; \colon \;   \left. v \right|_{\partial \mathcal{D} } =0     \right\} \label{def-H01}
\end{align}
\end{subequations}
 We define the following subspaces of $\left( H^1(\mathcal{D})\right)^2$:
\begin{subequations}
\begin{align}
\bfa{V}_{\bfzero} &:= \left\{ \bfu \in \left( H^{1}(\mathcal{D})\right)^2 \colon \; \bfu=\bfzero \quad  \mbox{on} \;\;\Gamma_D\right\}, \label{test_V0}\\
\bfa{V} &:= \left\{ \bfu \in \left( H^{1}(\mathcal{D})\right)^2 \colon \; \bfu=\bfu_0 \quad \mbox{on} \;\; \Gamma_D\right\}. \label{test_Vu0}
\end{align}
\end{subequations}

\subsection{Implicit theory of elasticity}
The primary objective of this paper is to develop a finite element model that examines the behavior of transversely isotropic, homogeneous elastic solids. Furthermore, it aims to present the crack-tip fields of the material body, the constitutive relationship of which is delineated through \textit{Rajagopal's theory of elasticity} \cite{rajagopal2003implicit,rajagopal2007elasticity,rajagopal2011conspectus,rajagopal2011non,rajagopal2014nonlinear}. In \cite{rajagopal2007elasticity}, Rajagopal extends the relationship established in Cauchy elasticity by recognizing that Cauchy stress $\bfT$ and the left Cauchy-Green stretch tensor $\bfB$ are related via an implicit equation.
 \begin{equation}\label{implicit_1}
\bfzero= \widetilde{\bfa{F}}(\bfB, \; \bfT).
\end{equation}
One can also consider the general subclass in the context of \textit{Rajagopal's theory of elasticity} \eqref{implicit_1}):
\begin{equation}\label{eq_B_T}
\bfB := \widehat{\bfa{F}}(\bfT). 
\end{equation} 
Use of the classical linearization assumption concerning small displacement gradients \cite{gurtin1982introduction}, the constitutive model \eqref{eq_B_T} yields a nonlinear \textit{strain-limiting} constitutive relationship through a corresponding response function $\bfa{F} \colon \Sym(\mathbb{R}^{2 \times 2}) \mapsto \Sym(\mathbb{R}^{2 \times 2})$,
\begin{equation}\label{SL2}
\bfeps =  \bfa{F}( \bfT), \quad \mbox{with} \quad \max_{\bfT \in \Sym} \| \bfa{F}(\bfT) \| \leq \lambda, \quad  \lambda>0.
\end{equation}
Following the earlier works \cite{yoon2021quasi,yoon2022CNSNS,lee2022finite,itou2017contacting,itou2017nonlinear,itou2018states}, we consider a special choice of the nonlinear response function $\bfa{F} \colon L^1(\Gamma_c; \, \Sym(\mathbb{R}^{2 \times 2})) \mapsto L^{\infty}(\Gamma_c; \, \Sym(\mathbb{R}^{2 \times 2}))$ of the form:
\begin{equation}\label{def-F}
\bfa{F}(\bfT) = \dfrac{ \mathbb{K}[\bfT]}{ \left( 1 + \beta^{{a}} { \| \mathbb{K}^{1/2}[\bfT] \|}^{\alpha}                   \right)^{1/\alpha}}, \quad \mbox{with} \quad {\sup_{\bfT \in \Sym} \| \bfa{F}(\bfT) \|  \leq  \dfrac{1}{\beta}},
\end{equation}
where $\beta \geq 0$ and $\alpha >0$ are the modeling parameters. The fourth-order compliance tensor $\mathbb{K}[\cdot]$ is the inverse of the elasticity tensor $\mathbb{E}[\cdot]$ which is defined as 
\begin{equation}
\mathbb{E}[\bfeps] := 2 \mu\bfeps + \lambda \, \tr(\bfeps) \, \bfI + \gamma \left( \bfeps \colon \bfa{M} \right) \, \bfa{M}
\end{equation}
where $\mu >0$ and $\lambda >0$ are Lam$\acute{e}$ coefficients, $\gamma >0$ and $\bfa{M}$ is a structural tensor which defines the orientation of the fibers in the solid \cite{Mallikarjunaiah2015,MalliPhD2015}.  Based on \cite{itou2017nonlinear,itou2018states}, the following properties of the tensor valued function $\bfa{F}$ are listed:
\begin{itemize}
\item[(i)] It is noted that $\| \bfa{F}(\bfT_1) \| \leq \lambda$ for all $\bfT_1 \in \Sym(\mathbb{R}^{2 \times 2})$. The tensor-valued function $\bfa{F}(\cdot)$ possesses a two-sided bound. This property implies the existence of a uniform bound on the strain values, specifically in the vicinity of concentrators such as crack tips, re-entrant corners, and v-notches.
\item[(ii)] $ \left( \bfa{F}(\bfT_1) - \bfa{F}(\bfT_2) \right) \colon  \left( \bfT_1 - \bfT_2 \right)  > 0$ for all $\bfT_1, \; \bfT_2 \in \Sym(\mathbb{R}^{2 \times 2})$. The tensor-valued function $\bfa{F}(\cdot)$ is said to be strictly monotone. 
\item[(iii)] ${\left( \bfa{F}(\bfT_1) - \bfa{F}(\bfT_2) \right) \colon  \left( \bfT_1 - \bfT_2 \right) } \leq \widehat{c}_1 \; \| \bfT_1 - \bfT_2  \|^2$, for all $\bfT_1, \; \bfT_2 \in \Sym(\mathbb{R}^{2 \times 2})$ and $\widehat{c}_1$ is a constant depends upon the modeling and material parameters. This property implies that the tensor-valued function $\bfa{F}(\cdot)$ is continuous. 
\item[(iv)] The function $\bfa{F}(\cdot)$ is  coercive, i.e. there exists constant $\widehat{c}_2$ such that $ {\left|  \bfa{v} \cdot \bfa{F}(\bfa{\Pi}) \bfa{v} \right| \geq \widehat{c}_2 \| \bfa{v} \|^2  }$ for all $\bfa{\Pi} \in \Sym(\mathbb{R}^{2 \times 2})$ and $\bfa{v} \; (\neq \bfa{0}) \; \in \mathbb{R}^{2}$. The constant $\widehat{c}_2$ depends upon model and material parameters and also on the dimension of the problem. 
\end{itemize}

The transversely isotropic models introduced in \eqref{SL2} are both invertible for sufficiently small values of $\beta$, and the inverted constitutive relationship is hyperelastic. The hyperelastic equivalent formulation of the aforementioned constitutive relationship \eqref{SL2} is  
\begin{equation}\label{invRel}
\bfa{T} := \Psi\left( \|  \mathbb{E}^{1/2}[\bfeps] \|       \right) \, \mathbb{E}[\bfeps], \quad \mbox{with} \quad \Psi(s) = \dfrac{1}{ \left(1 - (\beta \, s)^{\alpha}          \right)^{1/\alpha}}
\end{equation}

In the subsequent section, the relation \eqref{invRel} is employed to articulate the boundary value problems to investigate the crack-tip fields within a transversely isotropic strain-limiting elastic body. The objective of this study is to compare the predictions generated by the proposed model with those derived from the linearized elastic model (i.e., with $\beta =0$ in \eqref{SL2}) concerning the crack-tip fields, including stress, strain, and strain energy.

\section{Boundary value problem and existence of solution}\label{BVP}
Investigating cracks and fractures in transversely isotropic elastic materials is essential for various reasons, primarily owing to the extensive application of such materials in critical engineering contexts. This significance is accentuated by the prevalence of transversely isotropic properties in numerous natural and engineered materials. These materials include composites, wood, rocks, geological formations, and biological tissues. The presence of cracks compromises the integrity of structures and may culminate in catastrophic failures. Understanding crack behavior within these materials is imperative for designing safe and reliable structures. In sectors such as aerospace, where composite materials are prevalent, the prediction and prevention of crack propagation are crucial for ensuring aircraft safety. Similarly, comprehending the deterioration of concrete structures, which can be significantly influenced by cracking, is paramount in civil engineering. Therefore, an investigation of crack-tip fields in transversely isotropic materials is warranted. Subsequently, a well-posed boundary value problem is formulated, and a stable numerical algorithm is presented. 

Utilizing the inverted constitutive relationship and the balance of linear momentum, the following  boundary value problem is obtained: 
\begin{subequations}
\begin{align}
-\nabla \cdot \bfa{T}  &= \bff,  \quad \mbox{in} \quad \mathcal{D}, \quad   \mbox{with} \quad \bfeps =  \Bigg[\frac{\mathbb{K}[\bfa{T}]}{(1 + \beta^a { \|\mathbb{K}^{1/2}[\bfa{T}]\| }^{a})^{1/a}}\Bigg]  \label{qpde} \\
\bfT \bfn &= \bfa{g}, \   \mbox{on} \;\; \Gamma_N,  \label{NCond}\\
\bfu &= \bfu_0, \ \mbox{on} \;\; \Gamma_D, 
\end{align}
\end{subequations}
where $\bff \in  \left(L^{2}(\mathcal{D}) \right)^2$ is the body force term and $\bfa{g} \in  \left(  H^{3/2}(\Gamma_N) \right)^2$ is the traction. For the well-posedness of the above model, the following assumptions are made.
\begin{itemize}
\item[A1:]  The modeling parameters, $\beta$  and  $\alpha$, are presumed to be globally constant. Nevertheless, it is pertinent to verify through direct numerical simulations whether \(\beta \to 0^+\) allows the nonlinear model to recover the same predictions as the linear elastic constitutive relations.  The material parameters $\mu$ and $\lambda$ are also assumed to be constants. 
\item [A2:]  In the case when $\Gamma_D = \emptyset$, also the natural compatibility condition on the Neumann datum $\bfg$ and the source term $\bff$ yields
\begin{equation}
\bfzero = \int_{\mathcal{D}} \bff \, d\bfx + \int_{\partial \mathcal{D}} \bfg \, ds \quad \mbox{where} \quad  {\Gamma_D = \emptyset}. 
\end{equation}
\item[A3:]  The Dirichlet data $\bfu_0 \in \left( W^{1, \, 1}(\mathcal{D}) \right)^2$, with $\bfeps\left(\bfu_0(\bfx)\right)$ for almost every $\bfx \in \overline{\mathcal{D}}$ contained in a compact set in $\mathbb{R}^{2 \times 2}$. 
\end{itemize}

\begin{theorem} 
 Consider a bounded, connected, Lipschitz domain $\mathcal{D} \subset \mathbb{R}^2$ with a partitioned Lipschitz boundary. The boundary is composed of a relatively open Dirichlet boundary $\Gamma_D$ and a relatively open Neumann boundary $\Gamma_N$, such that $\Gamma_D \, \cap \, \Gamma_N = \emptyset$ and $\overline{\Gamma_D \cup \Gamma_N} = \partial \mathcal{D}$.  Given a vector field $\bff \colon \mathcal{D} \to \mathbb{R}^2$, a traction vector $\bfg: \Gamma_N \to \mathbb{R}^2$, a prescribed displacement $\bfu_0 \colon \Gamma_D \to \mathbb{R}^2$, and a bounded, nonlinear constitutive mapping $\bfa{F} \colon \Sym(\mathbb{R}^{2 \times 2}) \to \Sym(\mathbb{R}^{2 \times 2})$, determine the solution pair $(\bfu, \, \bfT)$, where $\bfu \colon \overline{\mathcal{D}} \to \mathbb{R}^{2}$ represents the displacement field and $\bfT \colon \overline{\mathcal{D}} \to \Sym(\mathbb{R}^{2 \times 2})$ represents the stress tensor, satisfying the following system of equations:
 \begin{align}
- \div \, \bfT &= \bff \quad \mbox{in} \quad \mathcal{D}, \notag \\
 \bfeps &= \bfa{F}(\bfT) := \dfrac{\bfT}{(1 + \left( \beta \, \| \bfT \| \right)^\alpha)^{1/\alpha}}, \; \alpha>0, \beta \geq 0 \quad \mbox{in} \quad \mathcal{D}, \notag \\
\bfu &= \bfu_0,  \;\; \mbox{on} \;\; \Gamma_D \\
\bfT \, \bfn &= \bfg,  \;\; \mbox{on} \;\; \Gamma_N \notag 
\end{align}
\noindent Weak Formulation and Existence: The existence of a solution pair $(\bfu, \, \bfT) \in  \left( W^{1, \; 1}(\mathcal{D})\right)^{2} \times \Sym({L}^1(\mathcal{D}))^{2 \times 2}$ satisfying the weak form is investigated:
\begin{equation}
\int_{\mathcal{D}} \bfT \cdot \bfeps(\bfw)  d\bfx = \int_{\mathcal{D}} \bff \cdot \bfw  d\bfx \quad \text{for all}  \quad \bfw \in  \left(C^{1}_{\Gamma_D}(\overline{\mathcal{D}})\right)^2.
\end{equation}
\noindent This weak formulation demonstrates the integral relationship between the stress tensor and the strain rate, establishing a basis for the existence of solutions within the specified function spaces.
\end{theorem}

The above theorem is examined in \cite{beck2017existence} and is congruent with the problem addressed in the present study, albeit with a minor variance in the definition of the Cauchy stress tensor as specified in \eqref{invRel}. The function $\Psi(\cdot)$ satisfies following lemma:
\begin{lemma}
Suppose $\mathbb{U} = \left\{ \omega \in \mathbb{R}^2 \colon  0 \leq |w| < c_1          \right\}$, and let $\omega_1, \; \omega_2 \in \mathbb{U}$, then there exists two positive constatns $c_2$ and $c_3$ such that the ofllowing inequalities hold:
\begin{align}
|\Psi(|\omega_1|)\,\omega_1-\Psi(|\omega_2|)\,\omega_2| & \leq c_2|\omega_1-\omega_2|, \\
(\Psi(|\omega_1|)\,\omega_1-\Psi(|\omega_2|)\,\omega_2, \omega_1-\omega_2)&\geq c_3 |\omega_1-\omega_2|^2.
\end{align}
Then the function $\Psi(\cdot)$ is both Lipschitz continuous and monotonic. 
\end{lemma}
Hence, the formulation above with the operator $\Psi(\cdot)$  is a uniform monotone operator exhibiting at most linear growth at infinity. Furthermore, we suppose that the materials' parameters $ \mu, \, \lambda$, modeling parameters $\beta, \, \alpha$, and the data vectors $\bfa{f}, \,\bfa{g}, \, \bfa{u}_0$ all conform to the stipulations established in the aforementioned theorem. Consequently, there exists a unique pair $(\bfu, \; \bfT) \in (W^{1,1}(\Omega))^2 \times \Sym({L}^1(\Omega)^{2 \times 2})$ pertaining to the continuous weak formulation delineated previously. Our proof for the designated selection of the strain tensor is derived from the corresponding argument presented in \cite{beck2017existence}.
\section{Galerkin-type finite element method}\label{fem}
\subsection{Continuous weak formulation }
This section proposes a finite element discretization of the BVP developed earlier. To obtain a well-posed weak formulation, we multiply the strong formulation \eqref{qpde} with the test function from $V_{\bfa{0}}$ (as in \eqref{test_V0}) and then integrate by parts using Green's formula together with the boundary conditions given in \eqref{NCond}, we arrive at the following weak formulation. 
\begin{cwf}
Given all the parameters, find $\bfu \in  \bfa{V}$, such that 
\begin{equation}
\label{eq:weak_formulation}
a(\bfu, \, \bfv) = L(\bfv), \quad  
\forall\, \bfv \in  \bfa{V}_{\bfa{0}},
\end{equation}
where, the bilinear term $a(\bfu, \, \bfv)$ and the linear term $L(\bfv)$ are defined as:
\begin{subequations}
\label{def:A-L}
\begin{align}
a(\bfu, \, \bfv) &= \int_{\mathcal{D}} \Psi( { \|\mathbb{E}^{1/2} [\bfeps(\bfu)  ]\|}) \, \mathbb{E}[\bfeps( \bfu) ] \colon \bfeps( \bfv) \;   d\bfx \, , \\
 L (\bfv) &= \int_{\mathcal{D}} \bff \cdot \, \bfv \; d\bfx \, + \int_{\Gamma_N} \bfa{g} \cdot \bfa{v} \, ds \, .
\end{align} 
\end{subequations} 
\end{cwf}
Assume that $\bff \in  \left(L^{2}(\mathcal{D}) \right)^2$, $\bfa{u}_0 \in (H^{1/2}(\Gamma_D))^2$, and $\bfa{g} \in (H^{3/2}(\Gamma_D))^2$, respectively. Then, the above problem \eqref{eq:weak_formulation} admits a unique weak solution \\
$\bfa{u} \in \mathcal{U}_s := \left\{  \bfa{v} \in  \left( H^2(\mathcal{D}) \cap W^{1, \, \infty}(\mathcal{D}) \right)^2 \colon \left. \bfa{v} \right|_{\Gamma_D} = \bfa{u}_0       \right\}$. The weak solution $\bfa{u}$ satisfies the following inequality
\begin{equation}
  \| \bfa{u}\|_{H^2(\Omega)}  \leq~ \widehat{c} \;  \big( \| \bff\|_{L^2(\mathcal{D})}+\| \bfa{u}_0 \|_{L^2(\Gamma_D)}+\| \bfa{g}\|_{L^2(\Gamma_N)}\big),\label{ctsubd}
\end{equation}
where $\widehat{c}$ denotes the regularity constant. 

\begin{remark}
Notably, the function  $\Psi(\cdot)$ exhibits nonlinearity; consequently, in our numerical simulations, we have employed Picard's iterative algorithm to address the nonlinearity. However, it is widely recognized that the convergence of Picard's method is contingent upon a "good initial" estimate. Therefore, our approach involves initially solving the linear problem (i.e., setting  $\beta = 0$), and subsequently utilizing the solution obtained as an initial estimate for the nonlinear Picard's iteration. 
\end{remark}

\subsection{Discrete finite element formulation}
Consider a collection of mesh partitions, denoted by $\left\{ \mathcal{T}_h \right\}_{h > 0}$, that dissect the closed domain $\overline{\mathcal{D}}$ into non-overlapping open subregions, labeled $\tau_i$ for $i$ ranging from $1$ to $N_h$. These subregions collectively cover the entire domain, meaning $\overline{\mathcal{D}}$ is the union of their closures: $\overline{\mathcal{D}}=\bigcup_{i=1}^{N_h} \bar{\tau}_i$. The size of each subregion $\mathcal{K}$ is characterized by its diameter, $h_{\mathcal{K}}$, and the overall mesh size, $h$, is defined as the maximum of these diameters across all subregions in $\mathcal{T}_h$, i.e., $h:= \max_{\mathcal{K} \in \mathcal{T}_h} h_\mathcal{K}$. Each mesh $\mathcal{T}_h$ consists of these mutually exclusive elements $\mathcal{K}$, and their union forms the entire domain: $\overline{\mathcal{D}} = \cup_{\mathcal{K} \in \mathcal{T}_h} \mathcal{K}$. Furthermore, the boundary edges/faces $\mathscr{E}_{bd,h}$ are then partitioned into two subsets: $\mathscr{E}_{D,h}$, representing Dirichlet boundary edges/faces, and $\mathscr{E}_{N,h}$, representing Neumann boundary edges/faces, such that their union yields the entire boundary set: $\mathscr{E}_{bd,h}=\mathscr{E}_{D,h} \cup \mathscr{E}_{N,h}$. Consequently, the complete set of all edges/faces, $\mathscr{E}_{h}$, is the combination of the interior and boundary edges/faces: $\mathscr{E}_{h}=\mathscr{E}_{int, h} \cup \mathcal{E}_{bd,h}$.

Define $\bfa{V}_{h}$ as the finite element space comprising piecewise continuous, bilinear, vector-valued functions over the mesh $\mathcal{T}_h$. Specifically, $\bfa{V}_{h}$ is given by:
\begin{equation}
\bfa{V}_{h}:= \left\{ \bfa{u}_h \in \left( C^{0}\left( \overline{\mathcal{D}} \right) \right)^2 \colon \left. \bfa{u}_h \right|_{\mathcal{K}} \in \left( \mathbb{Q}_{1}(\mathcal{K})\right)^2, \forall \mathcal{K} \in \mathcal{T}_h \right\} \subset \left(H^{1}(\mathcal{D})\right)^2.
\end{equation}
Here, $\mathbb{Q}_{1}(\mathcal{K})$ represents the space of bilinear polynomials defined on the element $\mathcal{K}$.  A  continuous-Galerkin semilinear form $a_h \colon \bfa{V}_{h} \times \bfa{V}_{h} \to \mathbb{R}$ and a linear functional $l_h \colon \bfa{V}_{h} \to \mathbb{R}$, are introdcued which together define the discrete weak formulation.
\begin{dwf}
To find $\bfa{u}_h \in \bfa{V}_{h}$ that satisfies the following equation for all test functions $\bfa{\varphi}_h \in \bfa{V}_{h}$ such that:
\begin{equation}\label{wf_1}
a_h(\bfa{u}_h; \bfa{u}_h, \, \bfa{\varphi}_h) = l_h(\bfa{\varphi}_h).
\end{equation}
The semilinear form $a_h(\bfa{u}_h; \bfa{u}_h, \, \bfa{\varphi}_h)$ and the linear form $l_h(\bfa{\varphi}_h)$ are defined as:
\begin{equation}\label{eqn:Jacobian}
a_h(\bfa{u}_h; \bfa{u}_h, \, \bfa{\varphi}_h) = \sum_{\tau_i \in \mathscr{T}_h}\int_{\tau_i} \Psi\left( \|\mathbb{E}^{1/2} [\bfeps(\bfu_h) ]\| \right) \, \mathbb{E}[\bfeps(\bfu_h)] \colon \bfeps( \bfa{\varphi}_h) \, d\bfa{x}
\end{equation}
and
\begin{equation}\label{eqn:RHS}
l_h(\bfa{\varphi}_h)= \sum_{\tau_i \in \mathscr{T}_h}\int_{\tau_i} \bff \cdot \bfa{\varphi}_h\,d\bfa{x} + \sum_{e_i \in \mathscr{E}_{N,h}} \int_{e_i} \bfa{g} \cdot \bfa{\varphi}_h\,ds.
\end{equation}
\end{dwf}
Subsequently,  two lemmas are examined that will yield Lipschitz continuity and strong monotonicity conditions for the semilinear form $a_h$ and the linear form $l_h$. 
\begin{lemma}
For $\bfa{u}_1, \; \bfa{u}_2 \in \bfa{V}_{h}$, the semi-linear form  $a_h$ is Lipschitz continuous (in the first argument), i. e. 
\begin{equation}
| a_h(\bfa{u}_1; \bfa{u}_1, \, \bfa{\varphi}_h)     -   a_h(\bfa{u}_2; \bfa{u}_2, \, \bfa{\varphi}_h)      | \leq k_1  \| \bfa{u}_1 - \bfa{u}_2 \| , \quad \mbox{for all} \quad \bfa{\varphi}_h \in \bfa{V}_{h}. 
\end{equation}
\end{lemma}
 
\begin{lemma}
For $\bfa{\varphi}_1, \; \bfa{\varphi}_2 \in \bfa{V}_{h}$, the semi-linear form  $a_h$ is strong monotonicity condition, i. e. 
\begin{equation}
| a_h(\bfa{\varphi}_1; \bfa{\varphi}_1, \, \bfa{\varphi}_1 - \bfa{\varphi}_2)     -   a_h(\bfa{\varphi}_2; \bfa{\varphi}_2, \, \bfa{\varphi}_1 - \bfa{\varphi}_2)      | \geq k_2  \| \bfa{\varphi}_1 - \bfa{\varphi}_2 \|.  
\end{equation}
\end{lemma}

At this juncture, the existence and uniqueness of the problem through the application of Riesz representation theory are examined.
\begin{theorem}
Assuming that the above two lemmas hold, then the discrete problem \eqref{wf_1} has a unique solution in  $\bfa{V}_{h}$. 
\end{theorem} 
The proof of the above theorem is a direct consequence of the Reisz representation theory \cite{manohar2024hp}. 

\section{Numerical results and discussion}\label{rd}
This section presents a numerical investigation into the behavior of "limiting" and "small" strains, which exhibit unbounded and escalating stresses near crack tips. A single crack in a transversely isotropic solid is examined to demonstrate the effectiveness and originality of our proposed modeling framework.  The central aim is to highlight the benefits of the nonlinear model in accurately capturing near-tip strain phenomena. A conventional, bilinear, continuous Galerkin finite element method is employed for the proposed direct numerical approximations, which is sufficient for the current objectives. Advanced discretization techniques and comprehensive \textit{a priori} error analysis will be the focus of future investigations; thus, one may employ the analysis developed in \cite{manohar2024hp}. The mathematical model uses the open-source finite element library \texttt{deal.II} \cite{2023dealii, dealii2019design}. All numerical results are obtained on structured meshes. The algorithm used for the nonlinear strain-limiting model computations is detailed in Algorithm~\ref{alg1}. To monitor the convergence of Picard's iterations, the residual  at each iteration is computed, and the residual is defined as
\begin{equation}\label{res}
R^n_{\bfa{u}^n} = \int_{\mathcal{D}}  \Psi\left( \|\mathbb{E}^{1/2} [\bfeps(\bfu^n_h) ]\| \right) \, \mathbb{E}[\bfeps(\bfu^n_h)] \colon \bfeps( \bfa{\varphi}_h) \, d\bfa{x}. 
\end{equation}

The overall algorithm to compute the numerical solution to the discrete weak formulation \eqref{wf_1} is given as follows: 

\begin{algorithm}[H]
\SetAlgoLined
\KwInput{ Mesh, parameters such as $\nu$, $\alpha$, $\beta$, iteration number $M$, initial guess $\bfa{u}^0$}
\KwOutput{Computed values of $\bfa{u}^M$, and postprocessed data}
 \While{$[\text{Iteration Number} < \text{Max. Number of Iterations}]$.AND.$[\text{Residual} > \text{Tol.}]$}{
  \If{\text{Iteration Number} ==0}{
For the initial guess, solve the linear problem by taking $\beta=0$ in \eqref{wf_1}
   }
  Assemble \eqref{wf_1} using the basis funcitons from $\bfa{V}_{h}$ \;
  Solve $\bfa{u}^{n}$ by knowing $\bfa{u}^{n-1}$  \;
  Calculate Residual (Equation~\eqref{res})\;
  \If{$\text{Residual}\leq\text{Tol.}$}{
   Break\;
   } }
 Save the converged solution for postprocessing \;
 \caption{Picard's Iterative algorithm to approximate the solution to nonlinear crack boundary value problem}
 \label{alg1}
\end{algorithm}
 
The primary goal of this research is to evaluate the accuracy and applicability of a proposed nonlinear mathematical model for transversely isotropic solids, comparing it to the traditional linearized elasticity model. This investigation seeks a more realistic description of material response, especially in the critical region surrounding crack tips. This leads to enhanced predictive models for crack growth under mechanical loads. The new response relations for the behavior of the transversely isotropic solid, represented by the vector-valued PDE system   \eqref{qpde}, are solved numerically using Picard's iterative algorithm (Algorithm \ref{alg1}). To ensure convergence, we utilize a tolerance of $10^{-6}$ and limit the iterations to a maximum of  $10$. 

\begin{figure}[H]
  \centering
   \includegraphics[width=0.5\linewidth]{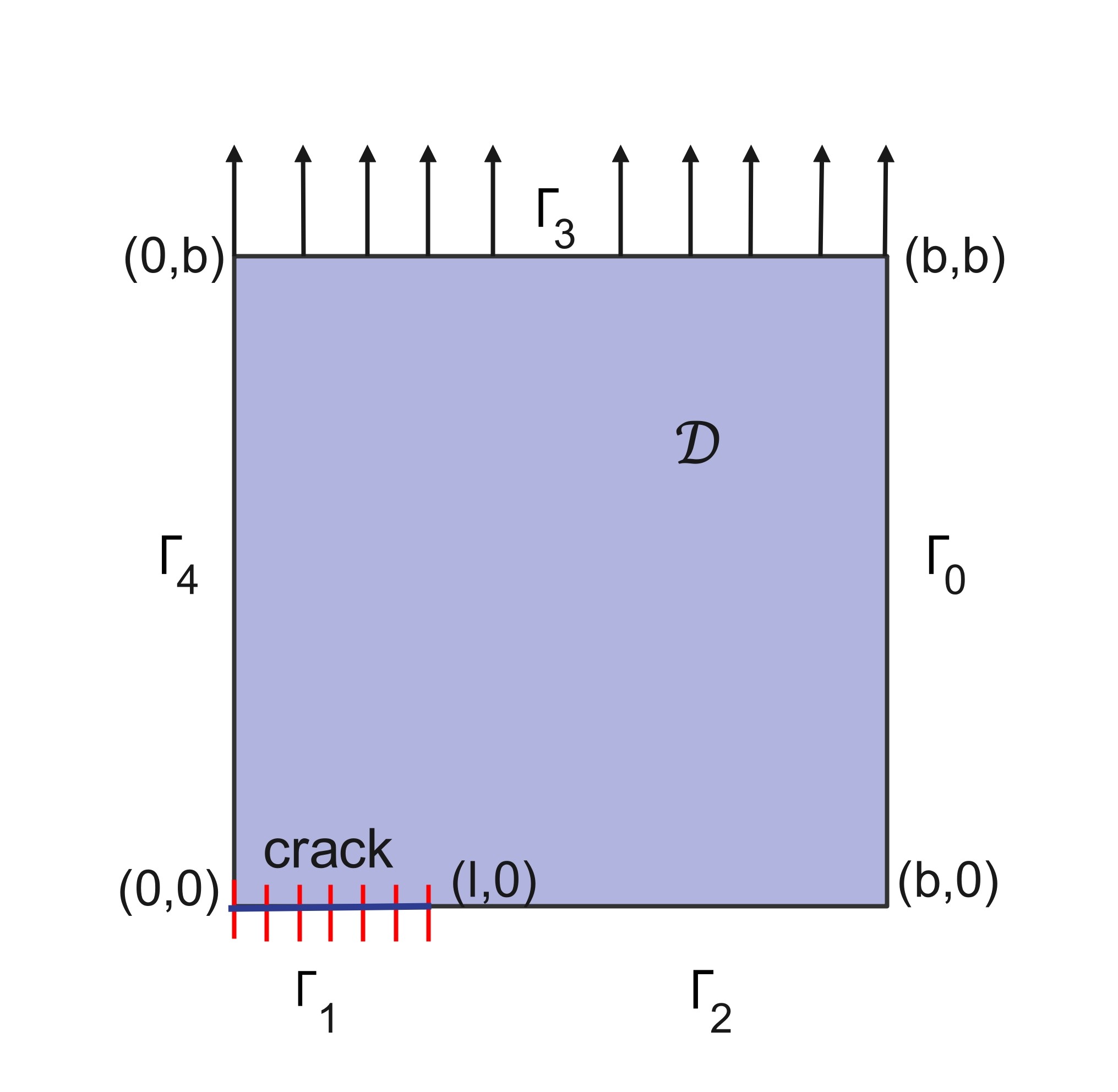}
        \caption{Computational domain}
   \label{fig:cd}
 \end{figure}

The numerical simulations utilized the computational domain depicted in Figure~\ref{fig:cd}. A mode-I (tensile) crack is modeled along the $x$-axis, specifically  $0 \leq x \leq 1, \; y=0$. A tensile load was applied to the upper boundary of the domain, and the right boundary is traction-free, while $u_1=0$ is imposed on the left boundary. The lower boundary was subjected to a vertical displacement constraint, setting $u_2=0$.
The material's transverse isotropy is characterized by examining two distinct fiber orientations, which are subsequently employed to modify the stress tensor. Notably, the symmetry axis of transversely isotropic materials, including fiber-reinforced composites, wood, and bone, often coincides with the orientation of the embedded fibers. This alignment is critical in ascertaining the material's overall response to applied loads.

\subsection{Fiber orientation along $x$-axis}
In this model, the material is treated as transversely isotropic, with the fibers oriented parallel to the x-axis, defining the axis of symmetry (parallel to the plane of the crack). For the fiber orientations, the structural tensor as $\bfa{M} = \bfa{e}_1 \bigotimes \bfa{e}_1$, with $\bfa{e}_1$ a unit vector along $x$-axis is considered. With a fixed Poisson's ratio, a Picard's iteration algorithm is executed. Table \ref{table1} presents the residual values at each iteration, demonstrating a significant decrease in the overall residual within the first six iterations.

\begin{table}[H]
  \centering
\begin{tabular}{|c|c|}
\hline
\textbf{Iteration No.} & $R^n_{\bfa{u}^n}$  \\
\hline
1 & 0.000798356 \\
\hline
2 & 0.000165949 \\
\hline
3 & 3.19741e-05 \\
\hline
4 & 5.64057e-06 \\
\hline
5 & 1.91659e-06 \\
\hline
6 & 1.3331e-06 \\
\hline
\end{tabular} 
  \caption{Residual computed at each iteration for the case of fiber's orientation is along the plane of the crack.}
  \label{table1}
\end{table}

In the remainder of this subsection, the numerical simulation results for fibers aligned parallel to the crack plane are presented. 

\begin{figure}[H]
    \centering
    \begin{subfigure}{0.3\linewidth}
        \centering
        \includegraphics[width=\linewidth]{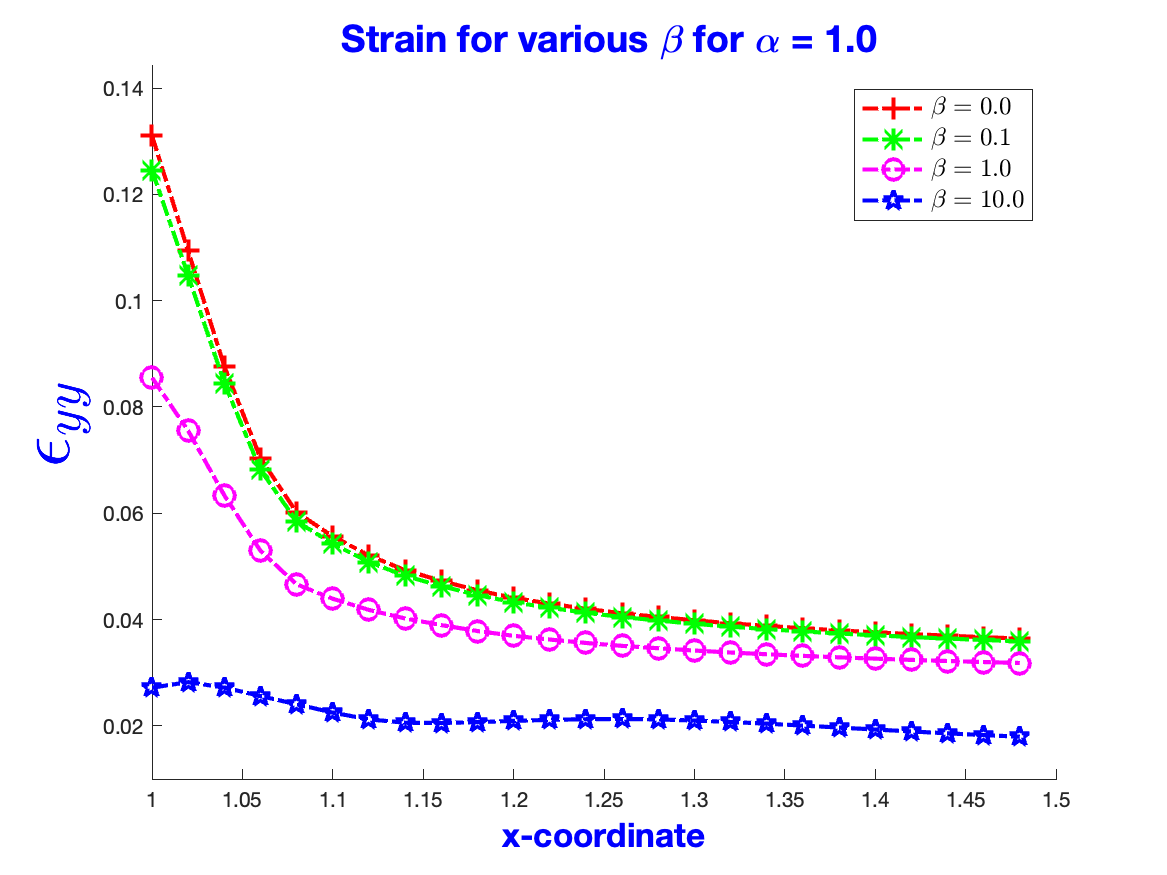}
        \caption{Strain for various $\beta$ for $\alpha = 1.0$ and $\sigma_{T} = 0.1$}
        \label{fig:strain_beta}
    \end{subfigure}
    \hfill
    \begin{subfigure}{0.3\linewidth}
        \centering
        \includegraphics[width=\linewidth]{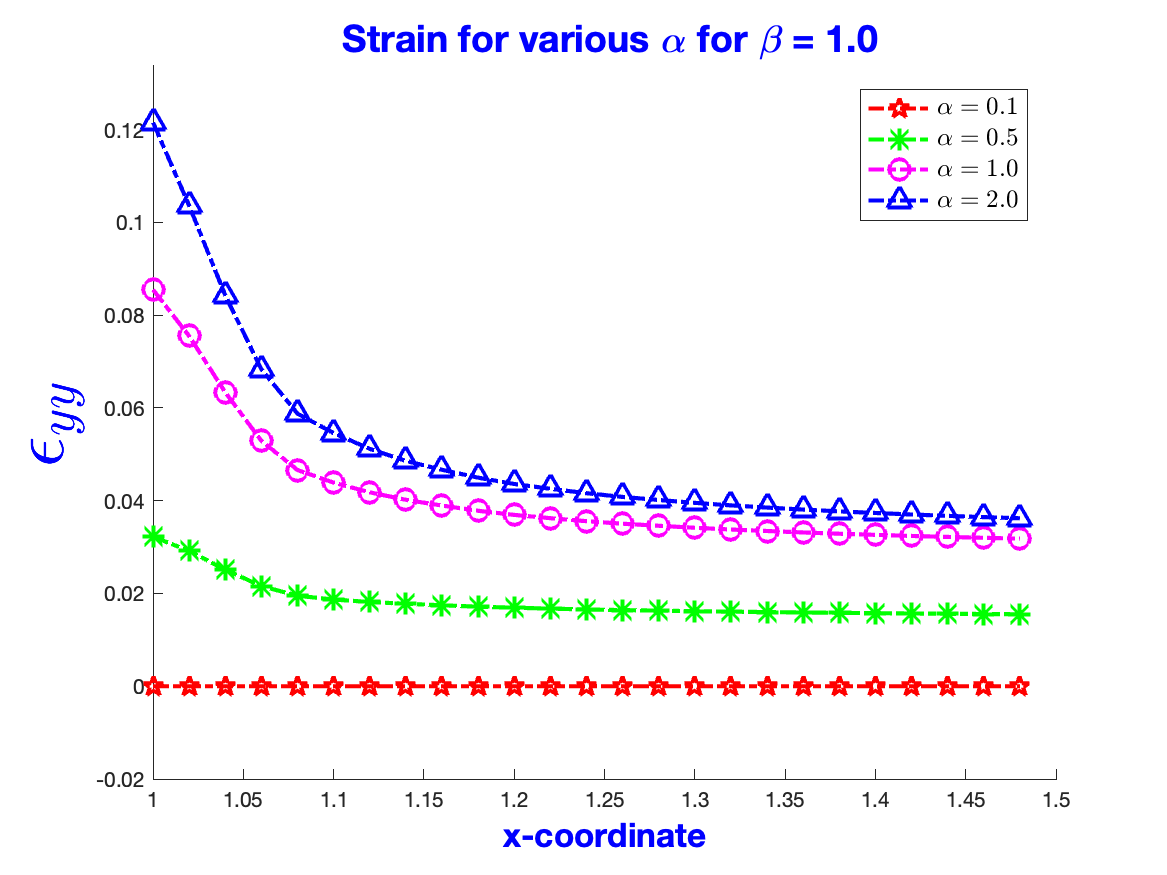}
        \caption{Strain for various $\alpha$ for $\beta = 1.0$ and $\sigma_{T} = 0.1$}
        \label{fig:strain_alpha}
    \end{subfigure}
    \hfill
    \begin{subfigure}{0.3\linewidth}
        \centering
        \includegraphics[width=\linewidth]{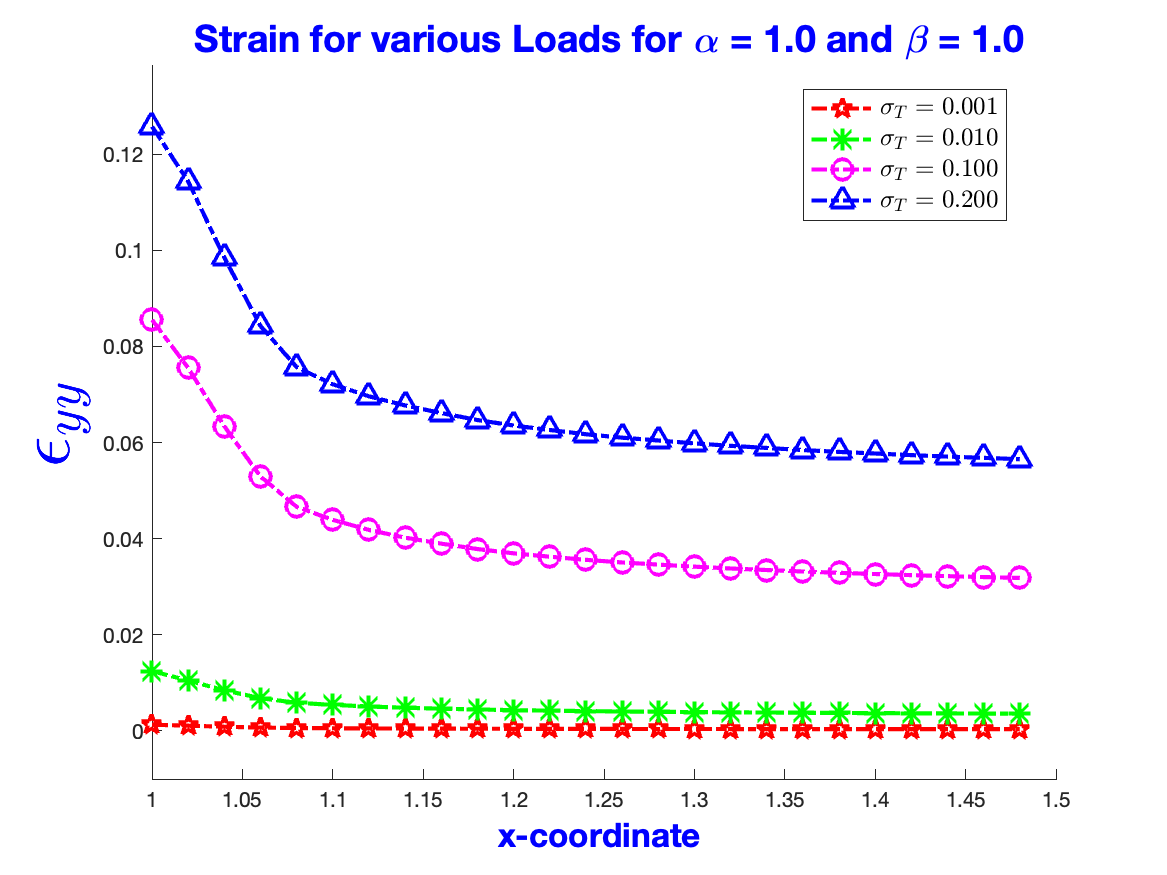}
        \caption{Strain for various $\sigma_{T}$ for $\beta = 1.0$ and $\alpha = 1.0$}
        \label{fig:strain_sigma}
    \end{subfigure}
    \caption{Strain plots for different parameter variations.}
    \label{fig:strain_combined}
\end{figure}

\begin{figure}[H]
    \centering
    \begin{subfigure}{0.3\linewidth}
        \centering
        \includegraphics[width=\linewidth]{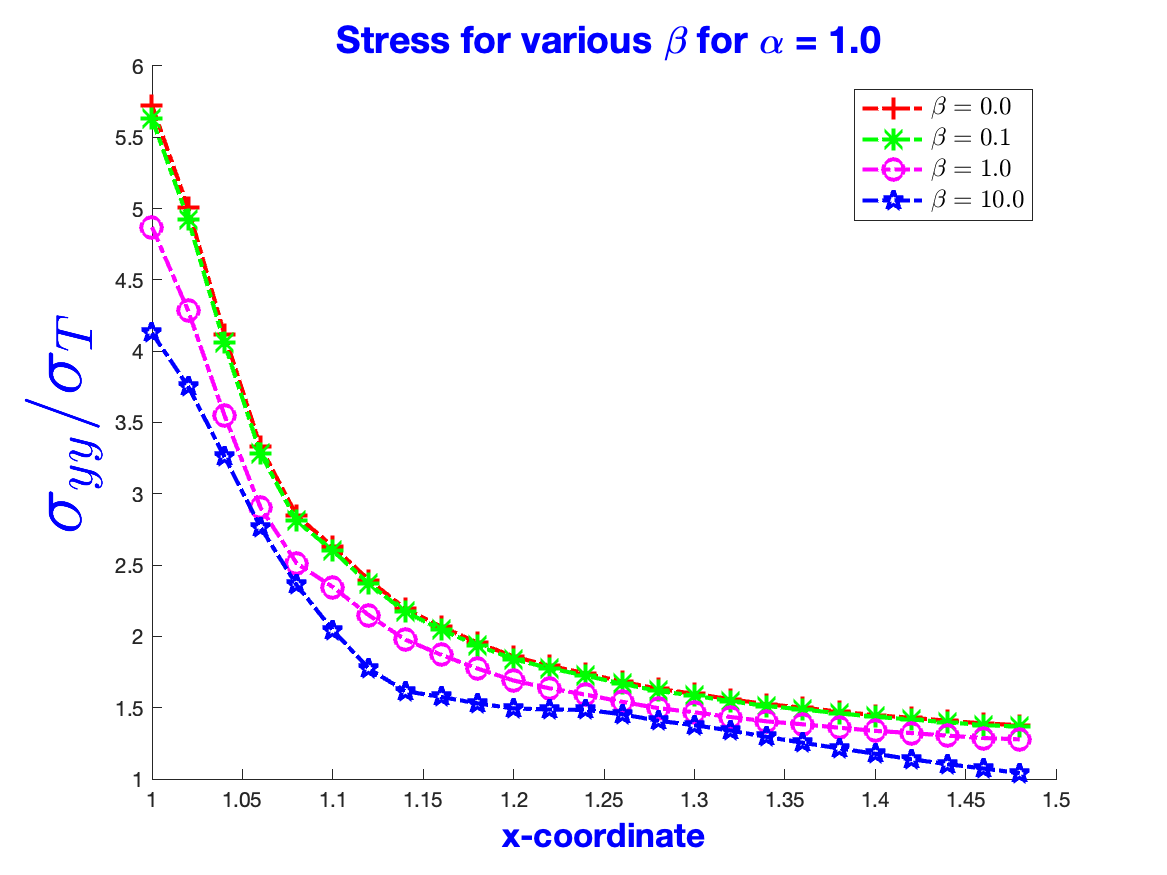}
        \caption{Stress for various $\beta$ for $\alpha = 1.0$ and $\sigma_{T} = 0.1$}
        \label{fig:stress_beta}
    \end{subfigure}
    \hfill
    \begin{subfigure}{0.3\linewidth}
        \centering
        \includegraphics[width=\linewidth]{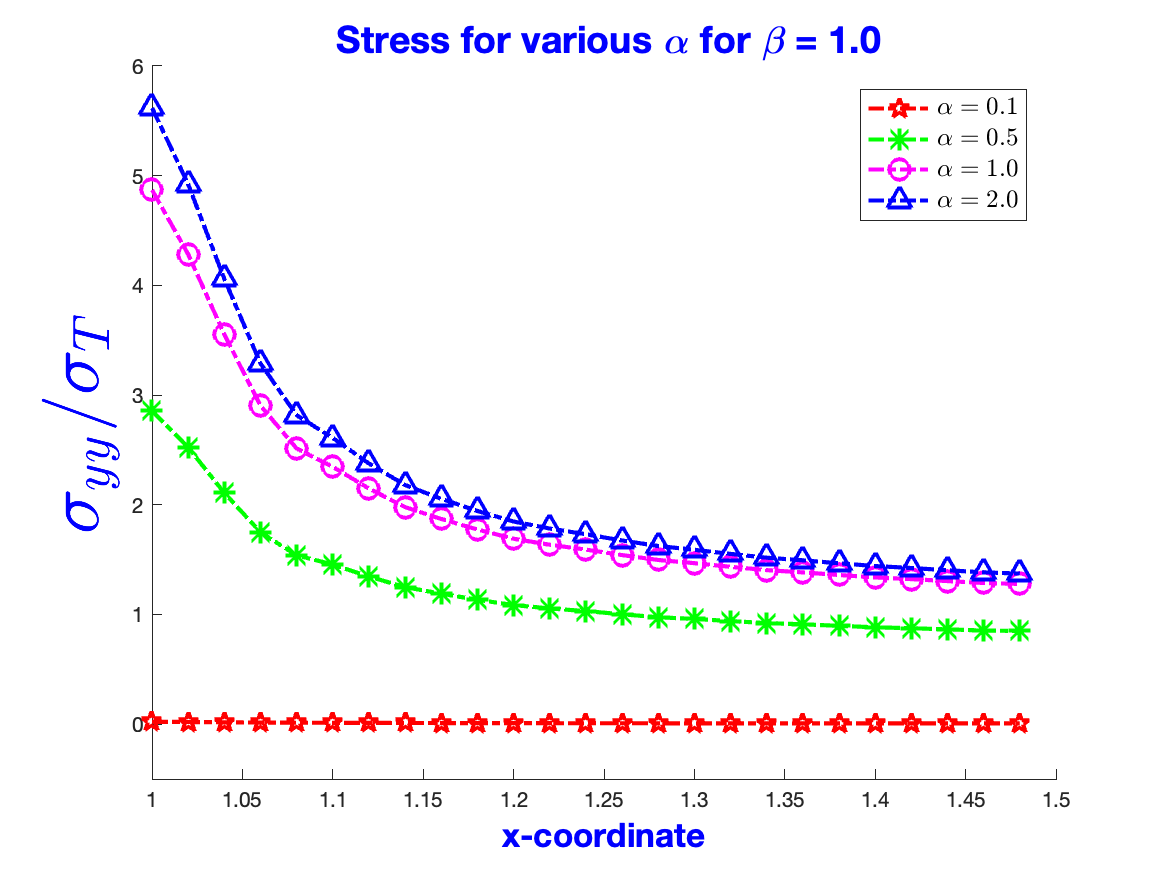}
        \caption{Stress for various $\alpha$ for $\beta = 1.0$ and $\sigma_{T} = 0.1$}
        \label{fig:stress_alpha}
    \end{subfigure}
    \hfill
    \begin{subfigure}{0.3\linewidth}
        \centering
        \includegraphics[width=\linewidth]{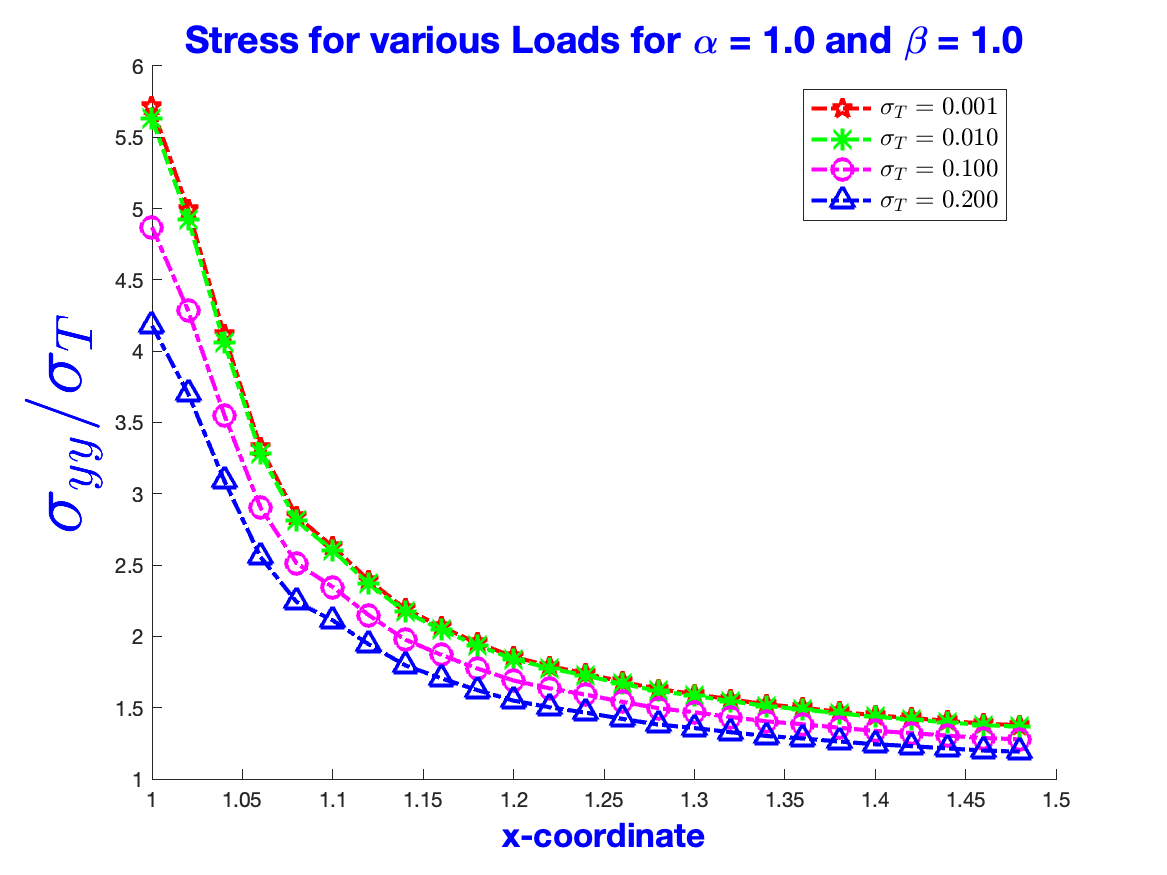}
        \caption{Stress for various $\sigma_{T}$ for $\beta = 1.0$ and $\alpha = 1.0$}
        \label{fig:stress_sigma}
    \end{subfigure}
    \caption{Stress plots for different parameter variations.}
    \label{fig:stress_combined}
\end{figure}

The nonlinear model's crack-tip strain and stress behaviors are presented in Figures~\ref{fig:strain_combined} and \ref{fig:stress_combined}. These figures demonstrate the sensitivity of both strain and stress to variations in $\beta, \; \alpha$, and $\sigma_T$, A key observation is the inverse relationship between $\beta$ and both crack-tip strain and stress. As $\beta$ increases, a significant reduction in these values is evident, confirming the strain-limiting characteristic of the proposed nonlinear formulation. While $\beta$ influences magnitude, substantial stress concentrations remain consistent across all investigated values, aligning with predictions from classical linearized models. Additionally, the figures reveal a direct proportional relationship between the applied top load 
$\sigma_T$ and the resulting crack-tip strain.

\begin{figure}[H]
    \centering
    \begin{subfigure}{0.3\linewidth}
        \centering
        \includegraphics[width=\linewidth]{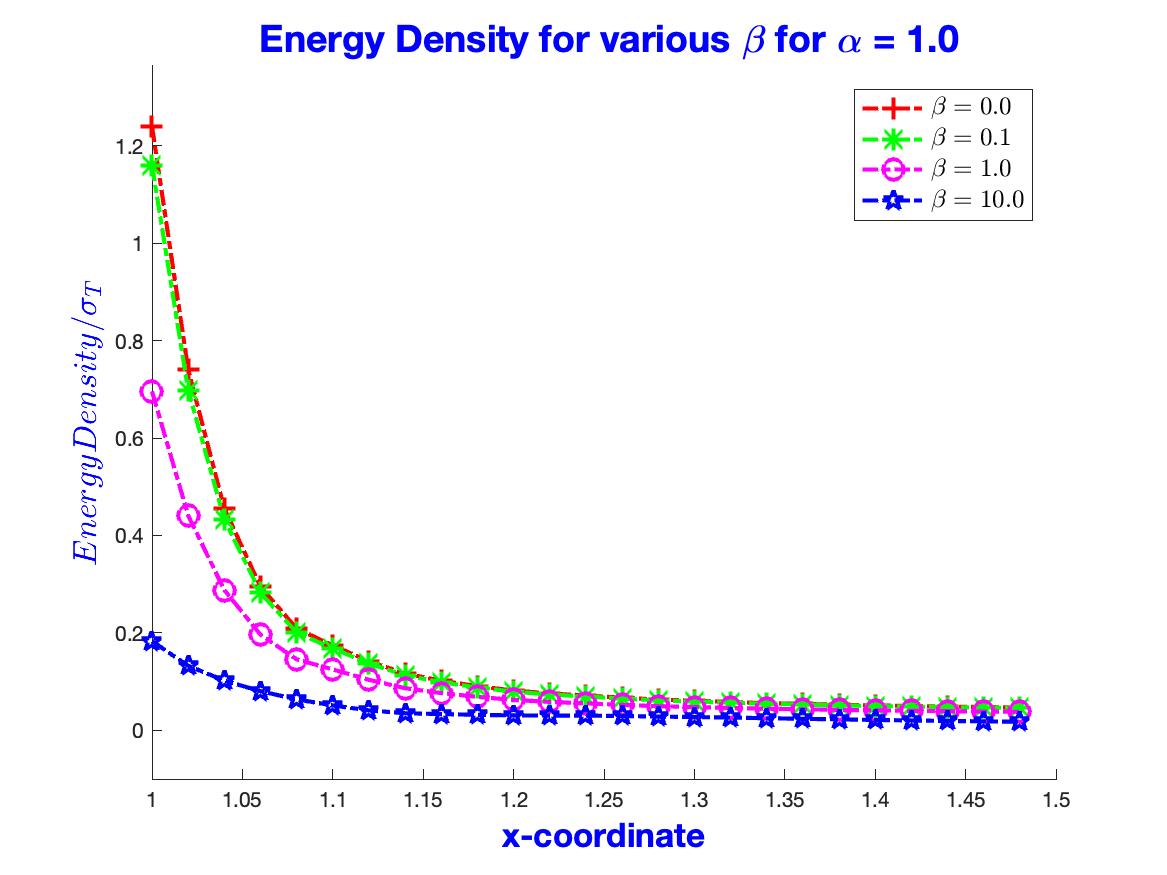}
        \caption{Energy density for various $\beta$ for $\alpha = 1.0$ and $\sigma_{T} = 0.1$}
        \label{fig:energy_density_beta}
    \end{subfigure}
    \hfill
    \begin{subfigure}{0.3\linewidth}
        \centering
        \includegraphics[width=\linewidth]{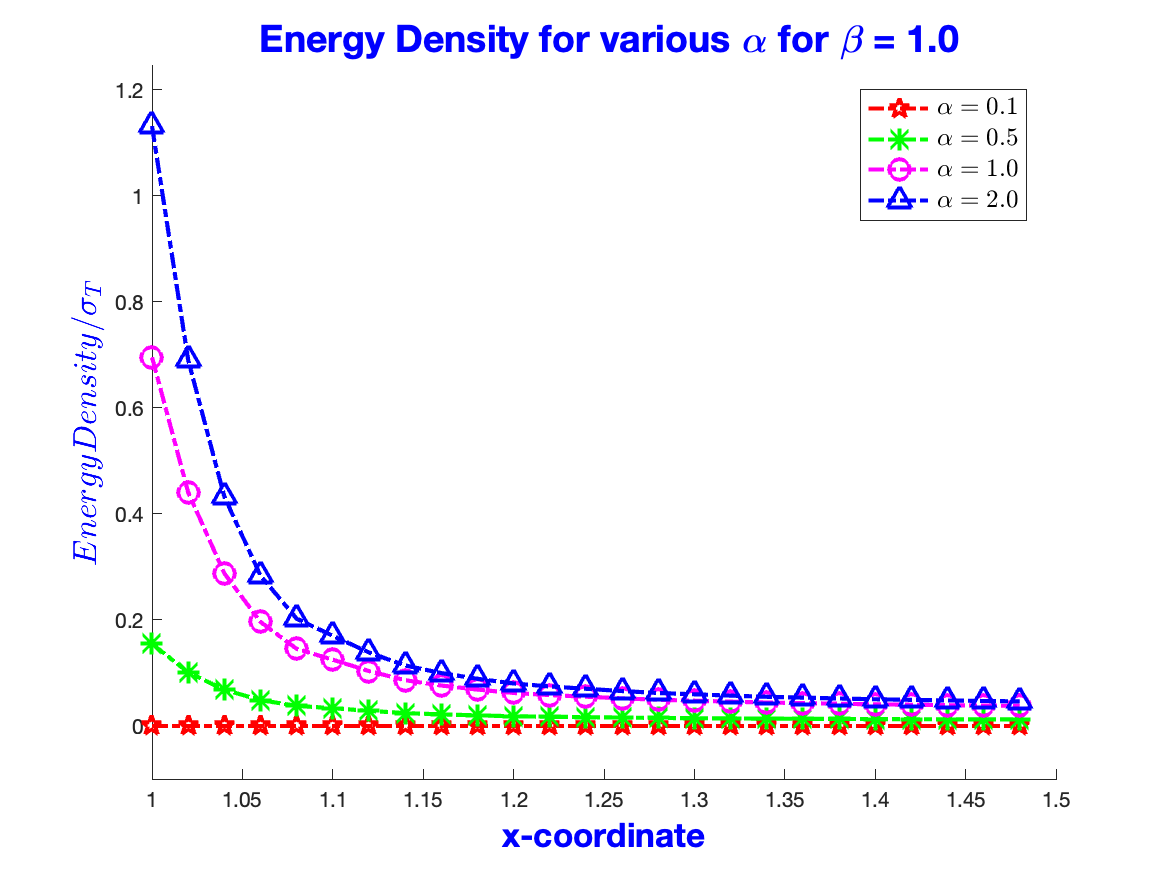}
        \caption{Energy density for various $\alpha$ for $\beta = 1.0$ and $\sigma_{T} = 0.1$}
        \label{fig:energy_density_alpha}
    \end{subfigure}
    \hfill
    \begin{subfigure}{0.3\linewidth}
        \centering
        \includegraphics[width=\linewidth]{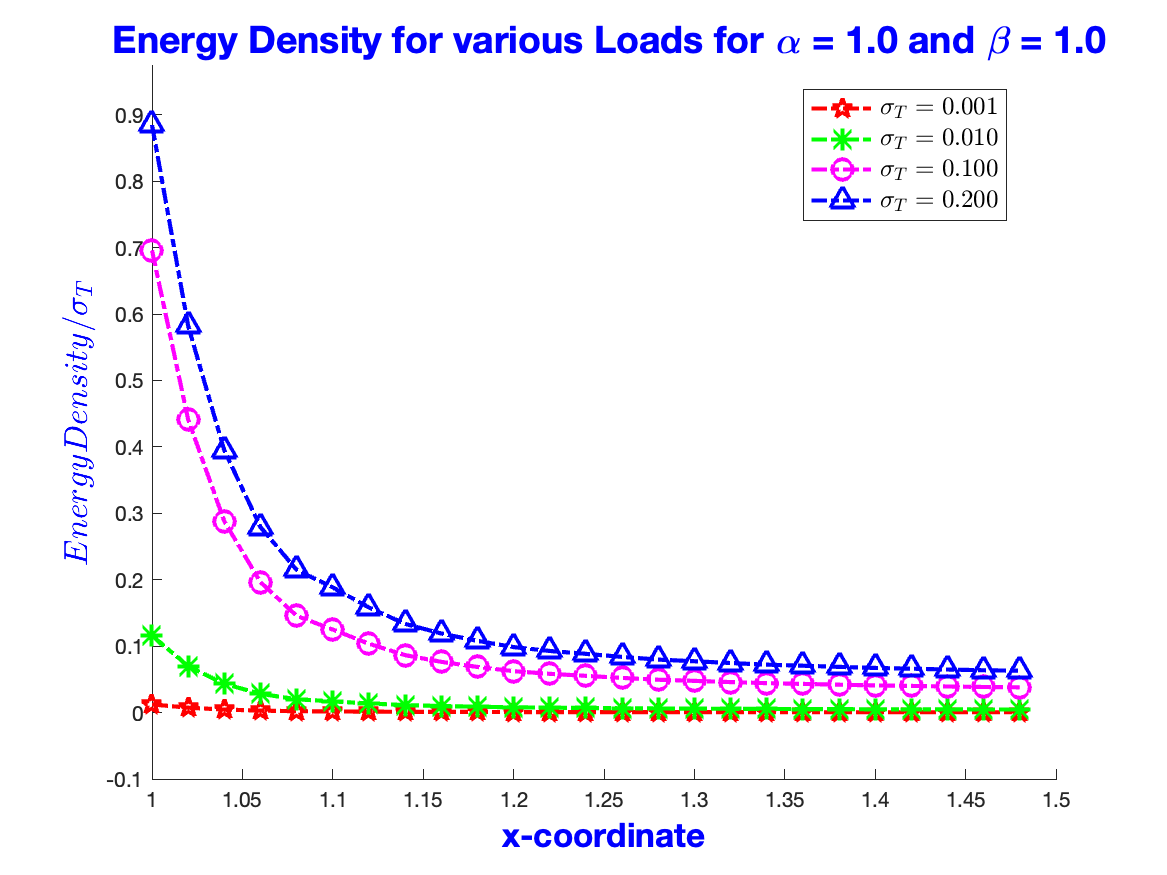}
        \caption{Energy density for various $\sigma_{T}$ for $\beta = 1.0$ and $\alpha = 1.0$}
        \label{fig:energy_density_sigma}
    \end{subfigure}
    \caption{Energy density plots for different parameter variations.}
    \label{fig:energy_density_combined}
\end{figure}

Figure~\ref{fig:energy_density_combined} depicts the strain-enrgy density ($\bfa{T} \colon \bfeps$) for various values of $\beta$, $\alpha$, and the top load $\sigma_T$. The strain-energy density at the crack tip is a critical parameter for predicting the onset of fracture and crack propagation. Analysis of the strain-energy density distribution near the crack tip provides insights into the local stress and deformation fields. It is observed that the strain-energy density decreases with increasing values of $\beta$, which means that the amount of energy stored within a material per unit volume is becoming less. However, an opposite effect is observed for the increasing values of $\alpha$. This observation indicates that a material demonstrates stress concentration without a corresponding strain concentration, suggesting a decoupling of stress and strain responses. This phenomenon implies that while the material undergoes localized high stresses, the associated deformation remains relatively uniform or constrained. Furthermore, it can be observed that geometric constraints within the material or structure may limit deformation, thereby preventing strain concentrations even in the presence of stress concentrations. It is worth noting that such results are obtained for a material model with geometric linearity.

\begin{figure}[H]
    \centering
    \begin{subfigure}{0.45\linewidth}
        \centering
        \includegraphics[width=\linewidth]{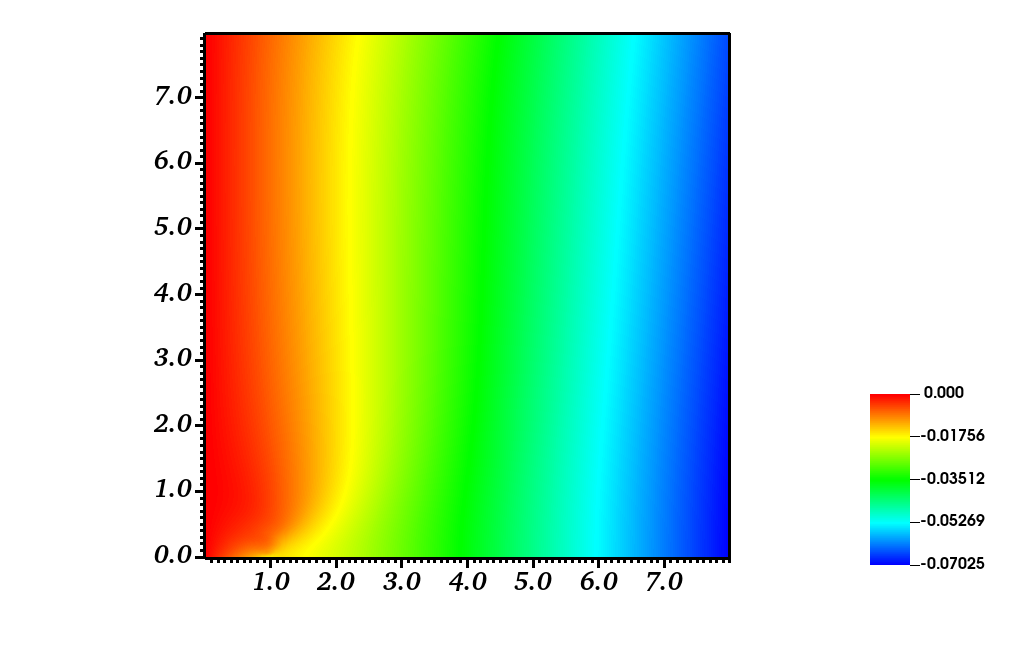}
        \caption{X-displacement for $\alpha = 1.0$, $\sigma_{T} = 0.1$, and $\beta = 10.0$}
        \label{fig:x_displacement}
    \end{subfigure}
    \hfill
    \begin{subfigure}{0.45\linewidth}
        \centering
        \includegraphics[width=\linewidth]{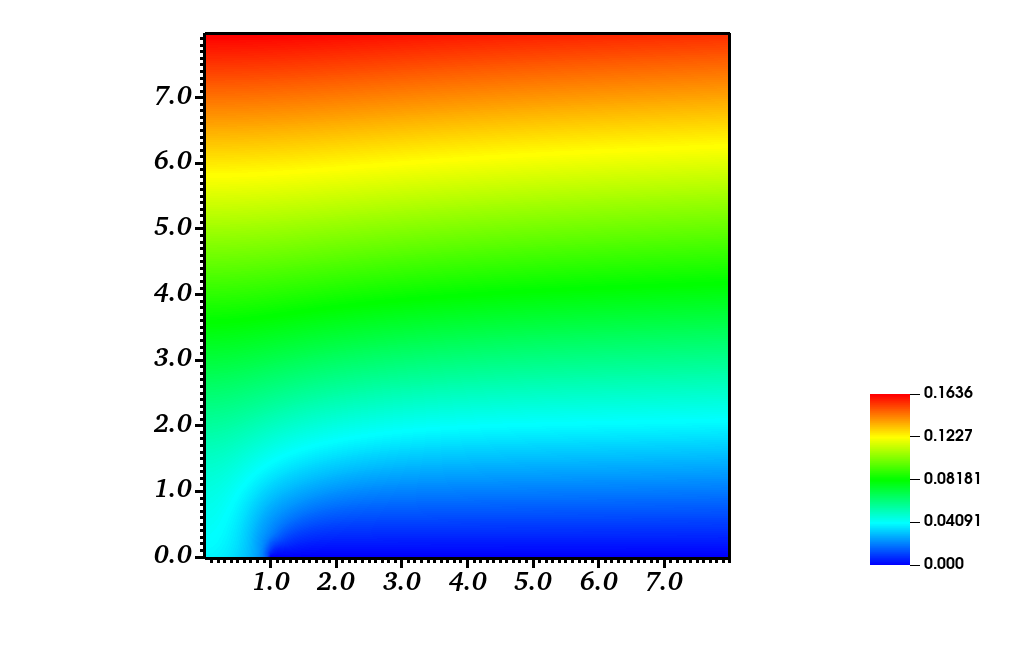}
        \caption{Y-displacement for $\alpha = 1.0$, $\sigma_{T} = 0.1$, and $\beta = 10.0$}
        \label{fig:y_displacement}
    \end{subfigure}
    
    \vspace{1em}
    
    \begin{subfigure}{0.45\linewidth}
        \centering
        \includegraphics[width=\linewidth]{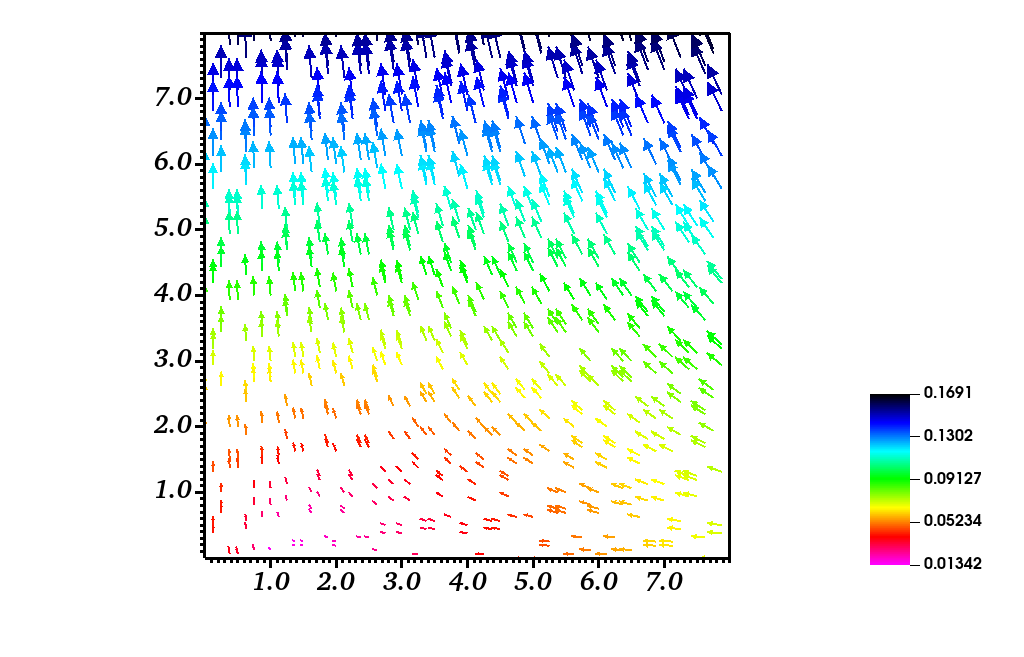}
        \caption{Vector-displacement for $\alpha = 1.0$, $\sigma_{T} = 0.1$, and $\beta = 10.0$}
        \label{fig:vector_displacement}
    \end{subfigure}
    \hfill
    \begin{subfigure}{0.45\linewidth}
        \centering
        \includegraphics[width=\linewidth]{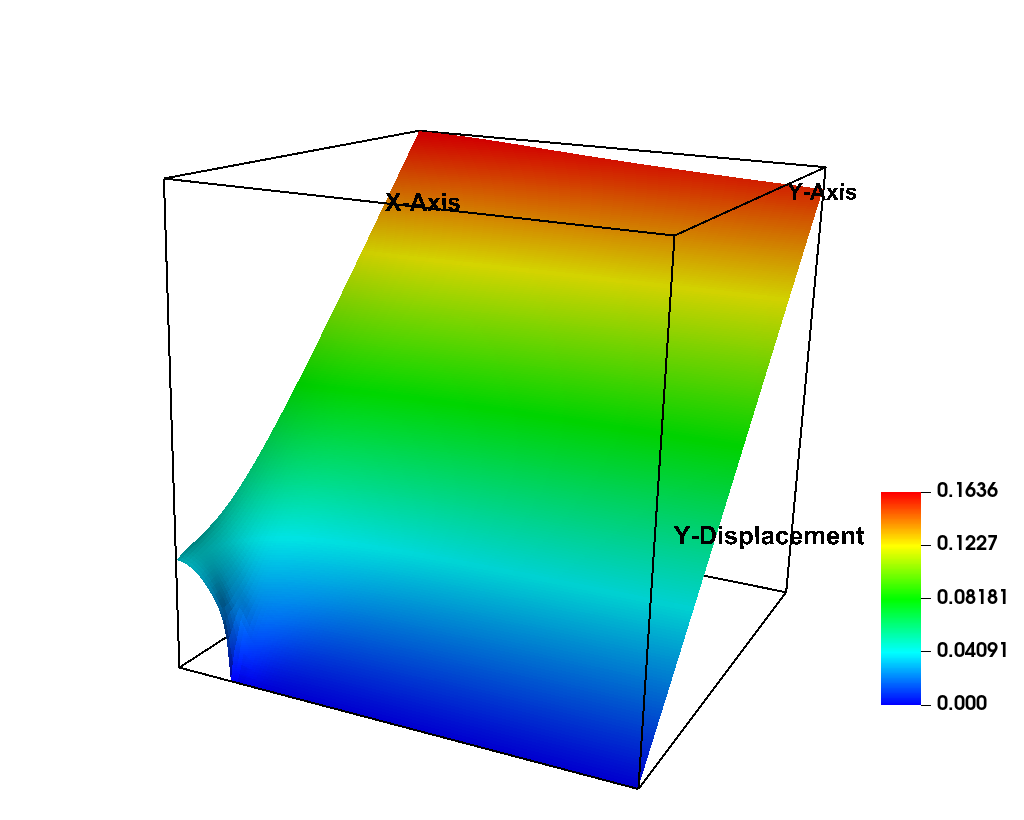}
        \caption{Y-displacement at $\alpha = 1.0$, $\sigma_{T} = 0.1$, and $\beta = 10.0$}
        \label{fig:y_displacement_3d1}
    \end{subfigure}
     \caption{Displacement plots for $\alpha = 1.0$, $\sigma_{T} = 0.1$, and $\beta = 10.0$ in various configurations.}
    \label{fig:displacement_combined1}
\end{figure}

Figure~\ref{fig:displacement_combined1} illustrates the graphical representation of the $x$ and $y$ displacements, the vector plot of the displacement $\bfa{u}$, and an elevated (three-dimensional) plot displaying the $y$-displacement, which also depicts the crack-opening profile. The crack-opening profile shown in \ref{fig:y_displacement_3d1} indicates that the crack opens with an elliptical profile with a blunt crack-tip. An elliptical opening signifies a gradual, curved deformation of the crack faces. This contrasts with a sharp, angular opening, which denotes a distinct stress distribution. Theoretically, a sharp crack tip would induce infinite stress concentration; however, this scenario is not physically realistic. A blunt crack tip mitigates the stress concentration, rendering the crack less susceptible to immediate propagation, which is commonly observed in materials with some degree of plasticity. The material yields at the crack tip, effectively blunting its profile.

\subsection{Fiber orientation perpendicular to $x$-axis}
Here, the material's behavior is modeled assuming transverse isotropy, with fibers aligned parallel to the $y$-axis (perpendicular to the $x$-axis), which served as the axis of symmetry and was coplanar with the crack.  The fiber orientation was represented by the structural tensor $\bfa{M} = \bfa{e}_2 \bigotimes \bfa{e}_2$, where $\bfa{e}_2$ is the unit vector parallel to the $y$-axis. A Picard's iterative scheme by utilizing a constant Poisson's ratio and variations in both the modeling parameters is implemented. The convergence of the algorithm is illustrated in Table \ref{table2}, which details a substantial reduction in the residual within the initial six iterations.

\begin{table}[H]
  \centering
\begin{tabular}{|c|c|}
\hline
\textbf{Iteration No.} &  $R^n_{\bfa{u}^n}$ \\
\hline
1 & 0.0006489 \\
\hline
2 & 0.000124101 \\
\hline
3 & 2.23515e-05 \\
\hline
4 & 3.61519e-06 \\
\hline
5 & 1.73245e-06 \\
\hline
6 & 1.419e-06 \\
\hline
\end{tabular}
 \caption{Values of the residual $R^n_{\bfa{u}^n}$  computed at each iteration for the current case of fiber's orientation perpendicular to the crack's plane.}
  \label{table2}
\end{table}

\begin{figure}[H]
    \centering
    \begin{subfigure}{0.3\linewidth}
        \centering
        \includegraphics[width=\linewidth]{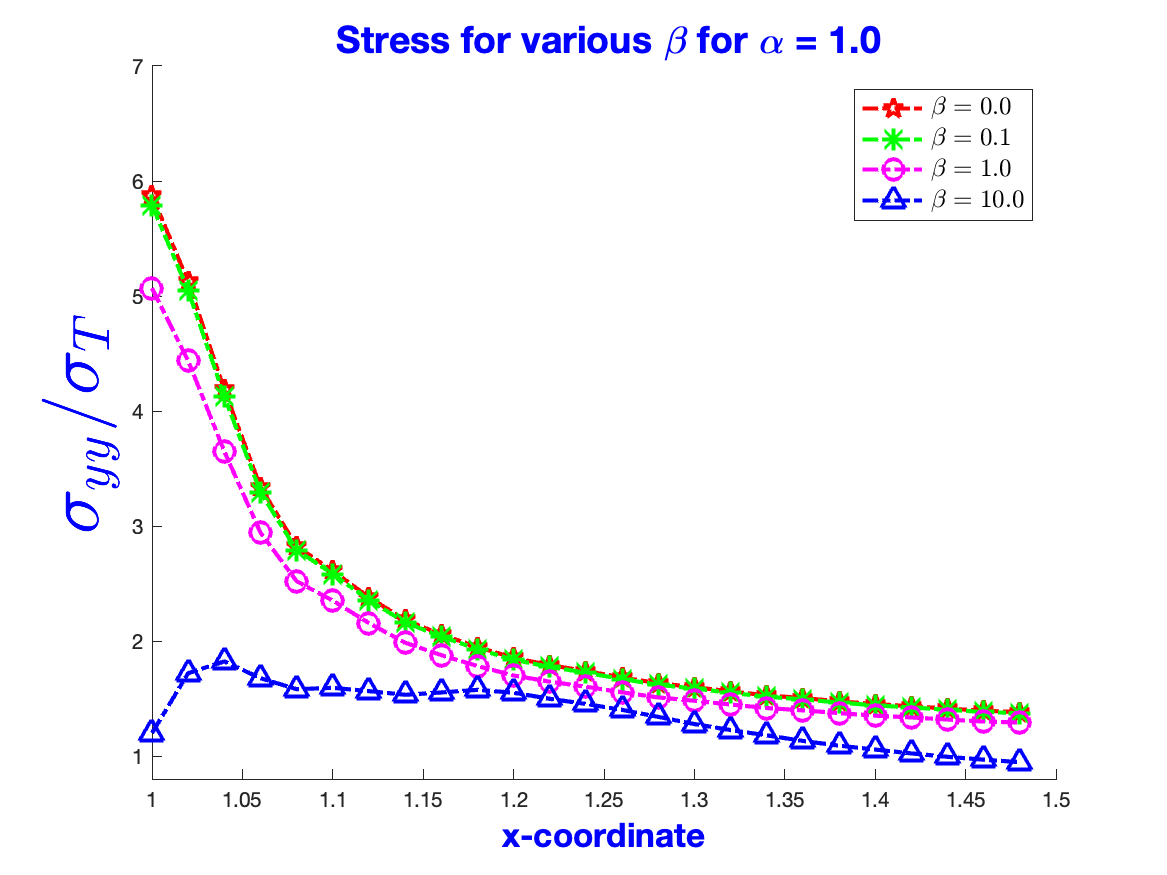}
        \caption{Stress for various $\beta$ for $\alpha = 1.0$ and $\sigma_{T} = 0.1$}
        \label{fig:stress_beta}
    \end{subfigure}
    \hfill
    \begin{subfigure}{0.3\linewidth}
        \centering
        \includegraphics[width=\linewidth]{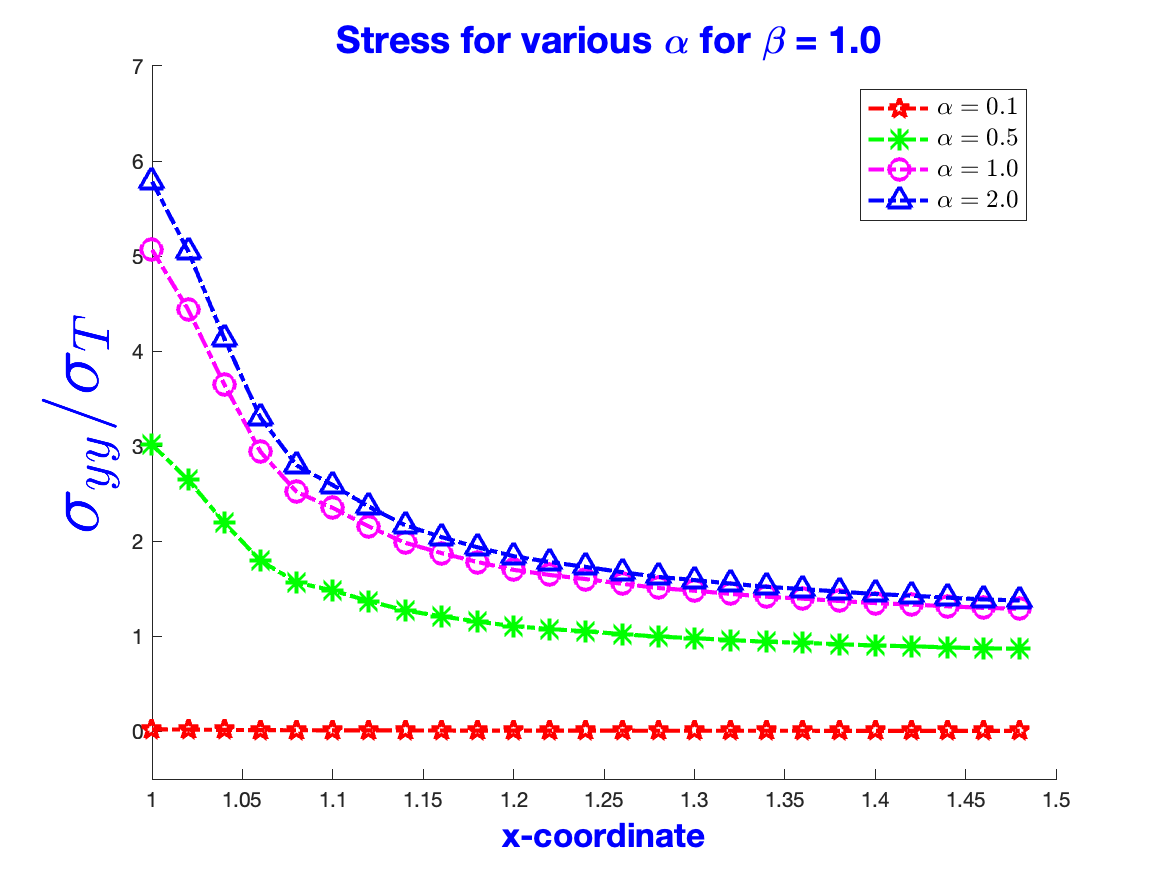}
        \caption{Stress for various $\alpha$ for $\beta = 1.0$ and $\sigma_{T} = 0.1$}
        \label{fig:stress_alpha}
    \end{subfigure}
    \hfill
    \begin{subfigure}{0.3\linewidth}
        \centering
        \includegraphics[width=\linewidth]{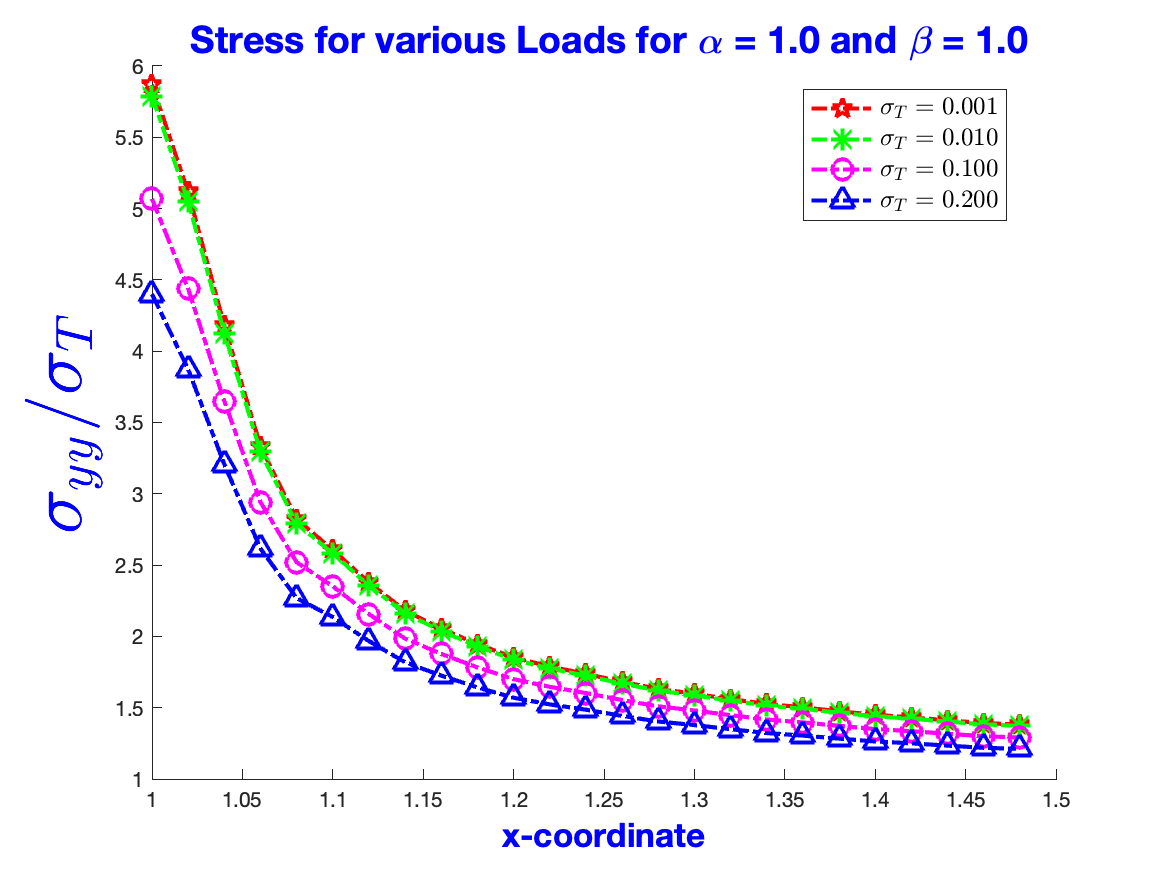}
        \caption{Stress for various $\sigma_{T}$ for $\beta = 1.0$ and $\alpha = 1.0$}
        \label{fig:stress_sigma}
    \end{subfigure}
    \caption{Stress plots for different parameter variations when the fiber directions are orthogonal to $x$-axis.}
    \label{fig:stress_combined_model2}
\end{figure}

\begin{figure}[H]
    \centering
    \begin{subfigure}{0.3\linewidth}
        \centering
        \includegraphics[width=\linewidth]{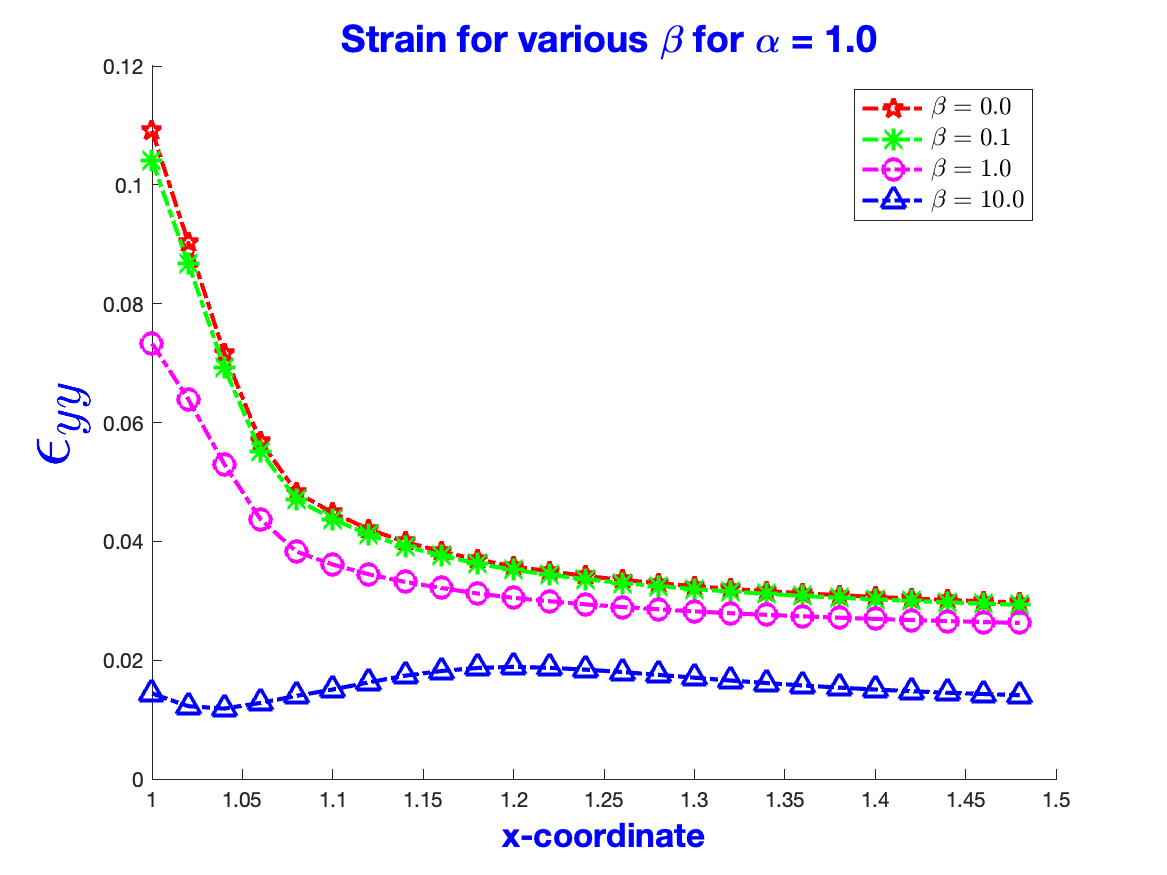}
        \caption{Strain for various $\beta$ for $\alpha = 1.0$ and $\sigma_{T} = 0.1$}
        \label{fig:strain_beta}
    \end{subfigure}
    \hfill
    \begin{subfigure}{0.3\linewidth}
        \centering
        \includegraphics[width=\linewidth]{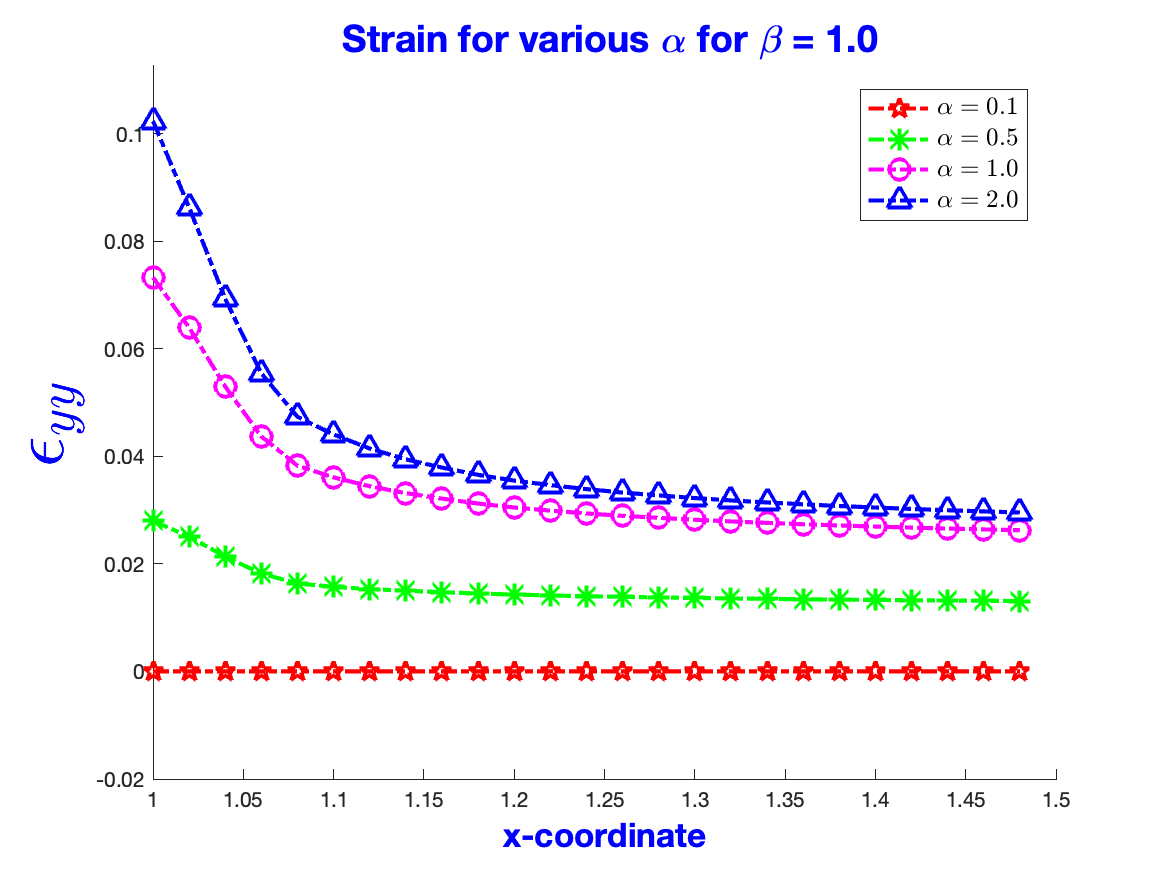}
        \caption{Strain for various $\alpha$ for $\beta = 1.0$ and $\sigma_{T} = 0.1$}
        \label{fig:strain_alpha}
    \end{subfigure}
    \hfill
    \begin{subfigure}{0.3\linewidth}
        \centering
        \includegraphics[width=\linewidth]{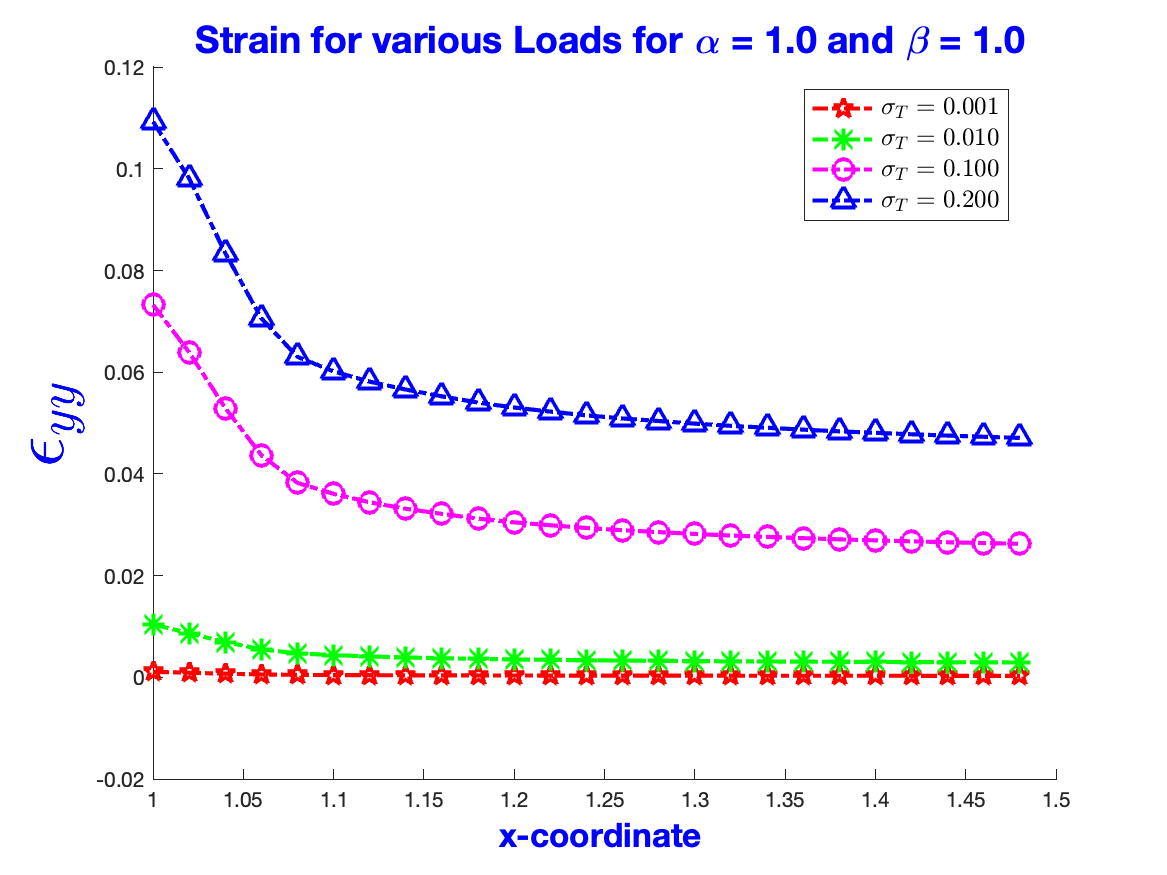}
        \caption{Strain for various $\sigma_{T}$ for $\beta = 1.0$ and $\alpha = 1.0$}
        \label{fig:strain_sigma}
    \end{subfigure}
    \caption{Strain plots for different parameter variations when the fiber directions are orthogonal to $x$-axis.}
    \label{fig:strain_combined_model2}
\end{figure}

Figures~\ref{fig:stress_combined_model2} and \ref{fig:strain_combined_model2} illustrate the stress and strain distributions along a radial path terminating at the crack tip, located at coordinates $(1, \, 0)$, for a material with fibers oriented orthogonally to the crack plane. Several intersting obeservations can be made from these figures, first, the crack-tip stresses behave pretty differently for inreasing values of the modeling parameter $\beta$, and for the value $\beta =10$ there is a sudden decrease in the crack-tip stress,  which indicates that the stress is being redistributed within the material. Lower crack-tip stresses suggest the crack is becoming more stable, meaning it's less likely to propagate.  However, increasing values of another modeling parameter $\alpha$ increases the neat-tip stress concentration, which favors the sudden evolution of the crack-tip. A slight decrease of the stresses in the neighborhood of the crack-tip for increasing top load, as depicted in the figure~\ref{fig:stress_sigma} suggests that more or less, the material behaves like strain hardening (although the strain hardening phenomenon is observed beyond the elastic limit to the material). 

The crack-tip strains decrease with the increasing values of the modeling parameter $\beta$ as shown in the firgure~\ref{fig:strain_combined_model2}, and this phenomenon occurs when a material or structure resists further deformation beyond a certain strain level. However, an opposite result is observed for the increasing values of another modeling parameter $\alpha$.  A slight increase trend is observed for the strains near the crack-tip for increasing top tensile load as shown in the figure~\ref{fig:strain_sigma}.

\begin{figure}[H]
    \centering
    \begin{subfigure}{0.3\linewidth}
        \centering
        \includegraphics[width=\linewidth]{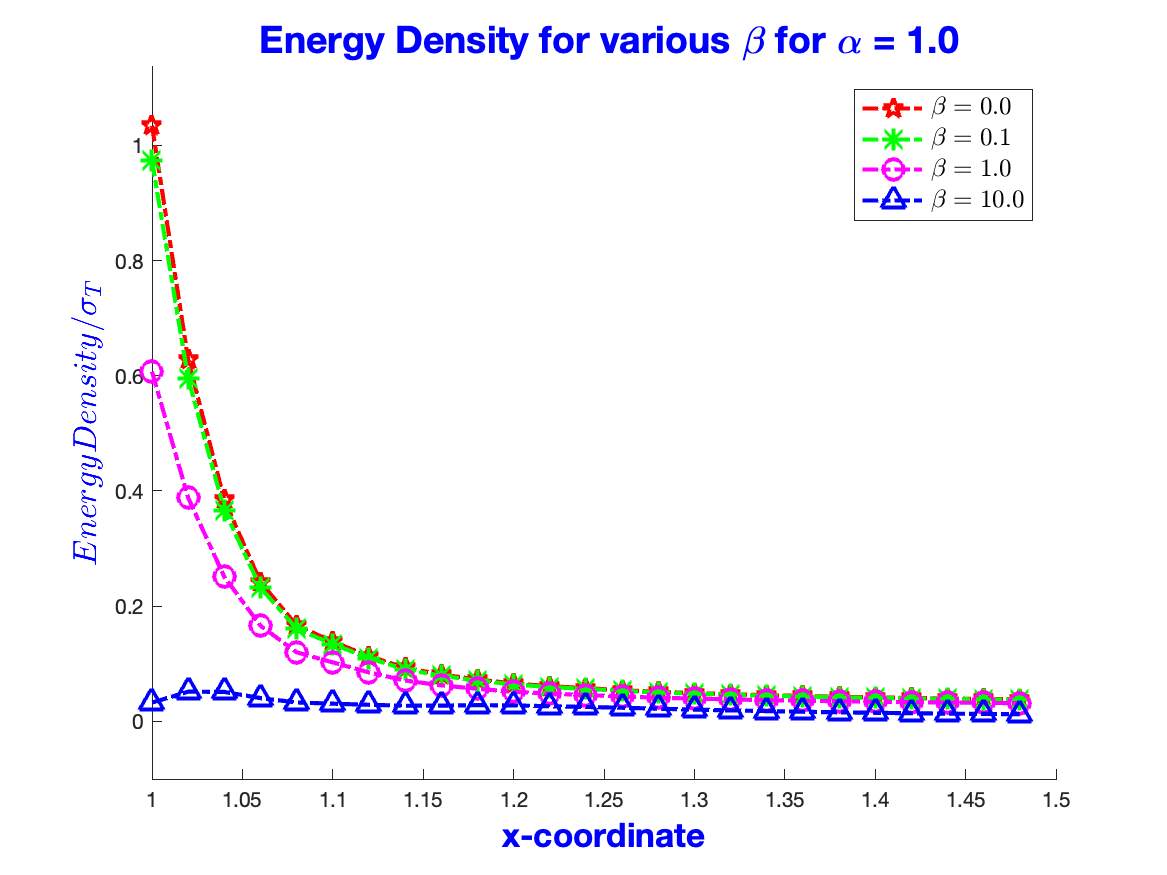}
        \caption{Energy density for various $\beta$ for $\alpha = 1.0$ and $\sigma_{T} = 0.1$}
        \label{fig:energy_density_beta}
    \end{subfigure}
    \hfill
    \begin{subfigure}{0.3\linewidth}
        \centering
        \includegraphics[width=\linewidth]{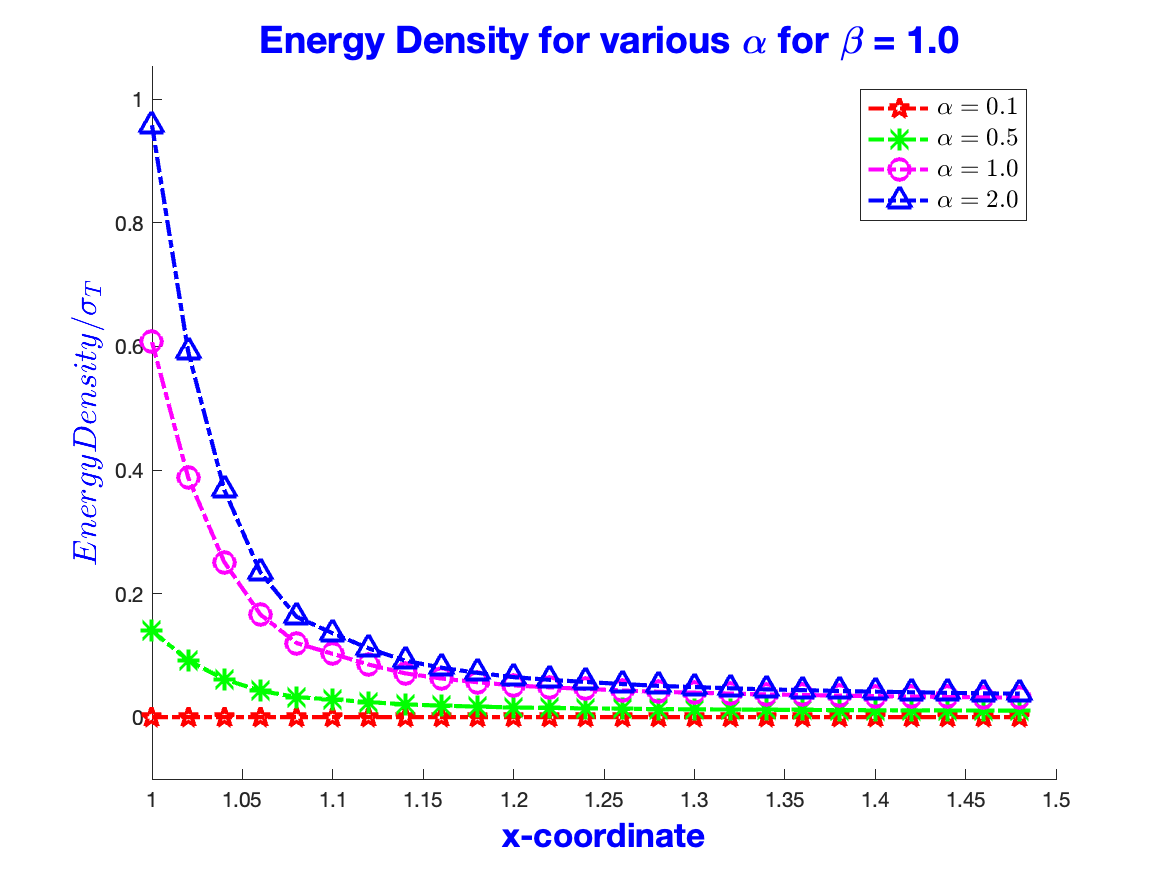}
        \caption{Energy density for various $\alpha$ for $\beta = 1.0$ and $\sigma_{T} = 0.1$}
        \label{fig:energy_density_alpha}
    \end{subfigure}
    \hfill
    \begin{subfigure}{0.3\linewidth}
        \centering
        \includegraphics[width=\linewidth]{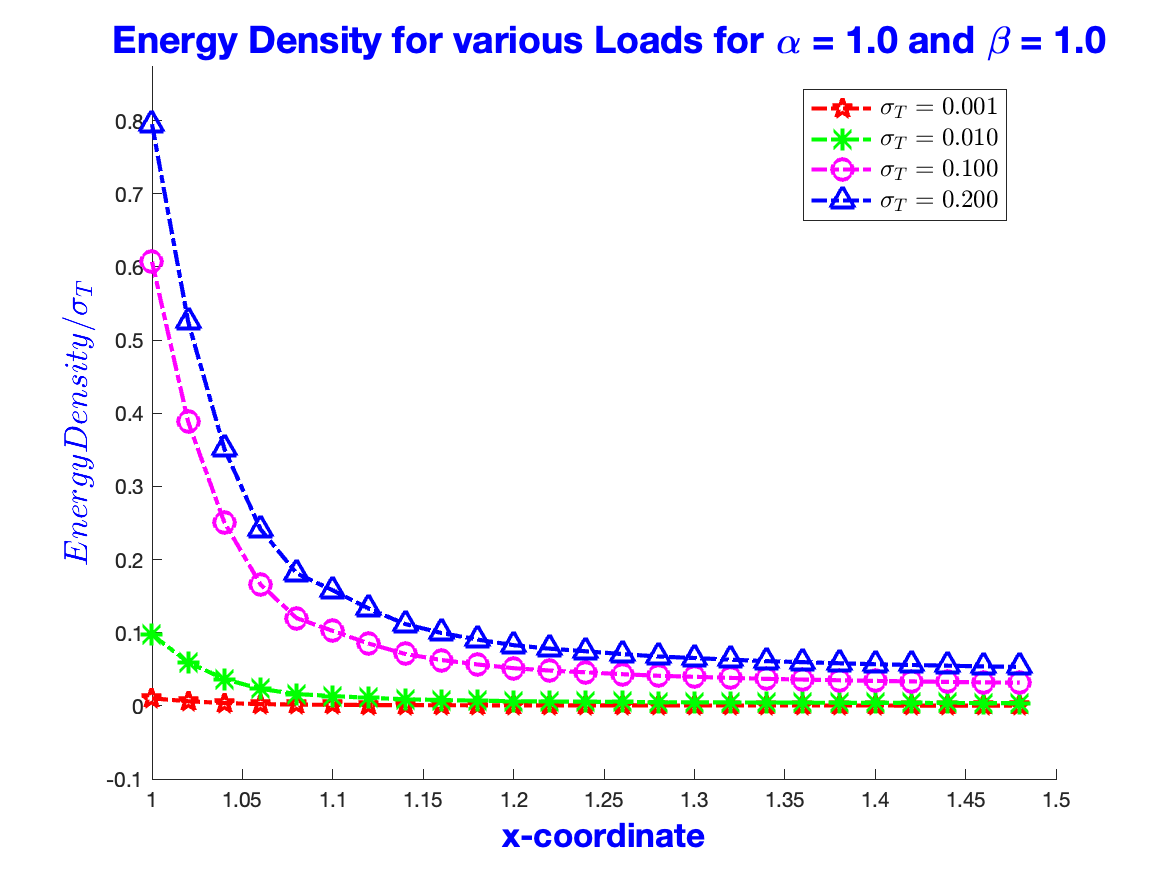}
        \caption{Energy density for various $\sigma_{T}$ for $\beta = 1.0$ and $\alpha = 1.0$}
        \label{fig:energy_density_sigma}
    \end{subfigure}
    \caption{Energy density plots for different parameter variations when the fiber directions are orthogonal to $x$-axis.}
    \label{fig:energy_density_combined_model2}
\end{figure}

Figure~\ref{fig:energy_density_combined_model2} illustrates the crack-tip energy density for various values of $\beta$, $\alpha$, and $\sigma_T$. The energy density exhibits a decreasing trend with increasing values of $\beta$, while an inverse behavior is observed for the increasing values of $\alpha$ and $\sigma_T$. An increase in strain energy density indicates that the material is undergoing more significant deformation, thereby reinforcing the crack-tip in accordance with the linearized theory of elastic fracture mechanics. Conversely, a decrease signifies that the material is releasing the potential energy it had previously accumulated from the deformation induced by the static loading conditions considered in this study. Typically, a reduction in strain energy density is associated with a decline in the internal stress within the material, which is evident from the figure showing the increase in $\beta$ values. 

\begin{figure}[H]
    \centering
    \begin{subfigure}{0.45\linewidth}
        \centering
        \includegraphics[width=\linewidth]{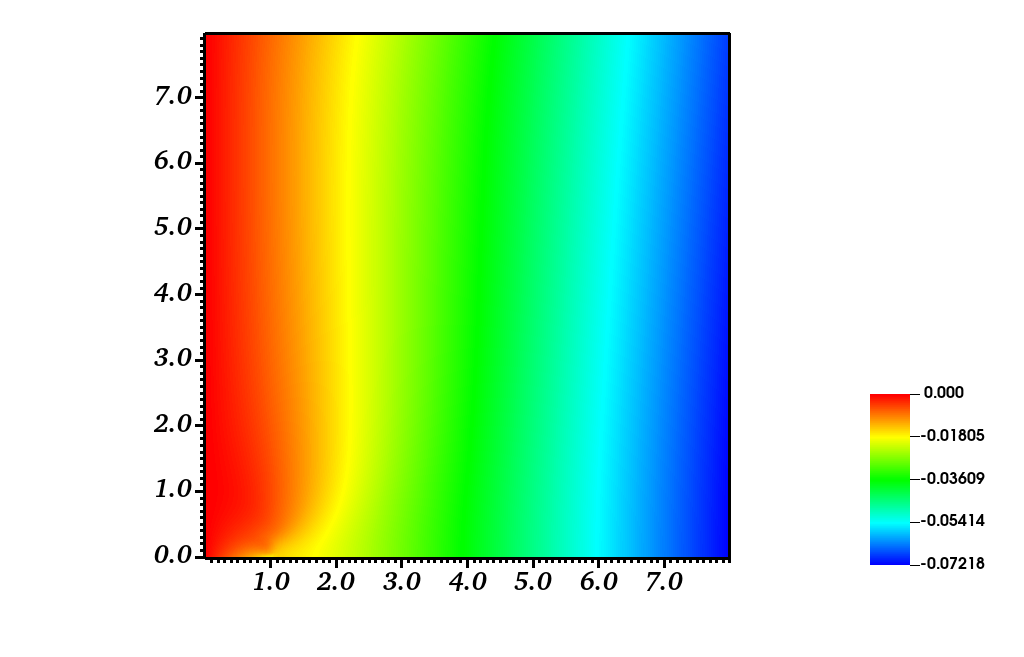}
        \caption{X-displacement at $\alpha = 1.0$, $\sigma_{T} = 0.1$, and $\beta = 10.0$}
        \label{fig:x_displacement}
    \end{subfigure}
    \hfill
    \begin{subfigure}{0.45\linewidth}
        \centering
        \includegraphics[width=\linewidth]{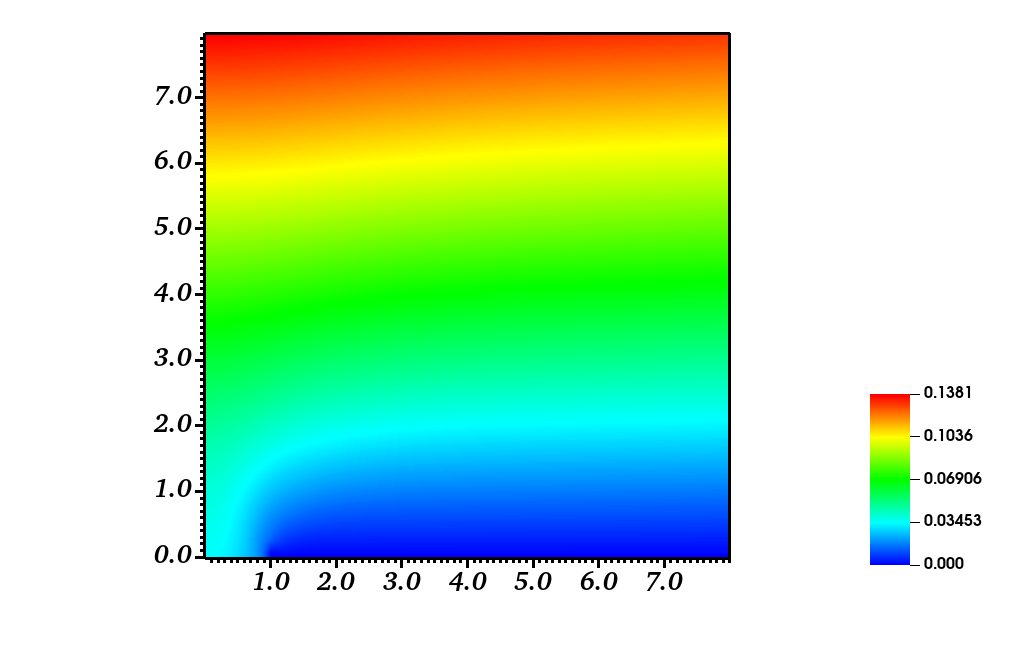}
        \caption{Y-displacement at $\alpha = 1.0$, $\sigma_{T} = 0.1$, and $\beta = 10.0$}
        \label{fig:y_displacement}
    \end{subfigure}
    
    \vspace{1em}
    
    \begin{subfigure}{0.45\linewidth}
        \centering
        \includegraphics[width=\linewidth]{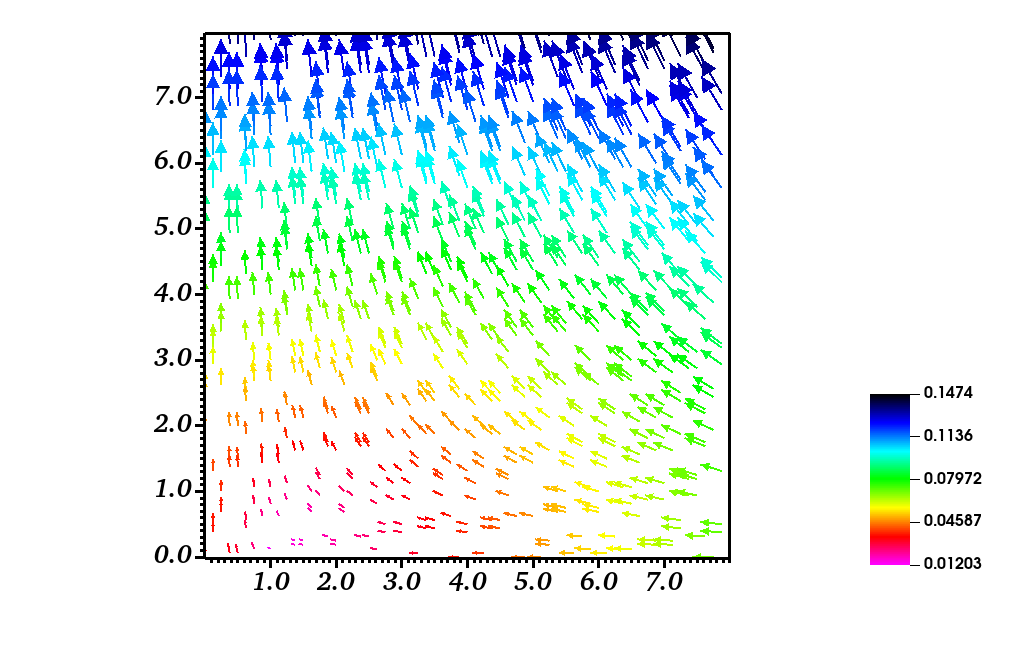}
        \caption{Vector-displacement at $\alpha = 1.0$, $\sigma_{T} = 0.1$, and $\beta = 10.0$}
        \label{fig:vector_displacement}
    \end{subfigure}
    \hfill
    \begin{subfigure}{0.45\linewidth}
        \centering
        \includegraphics[width=\linewidth]{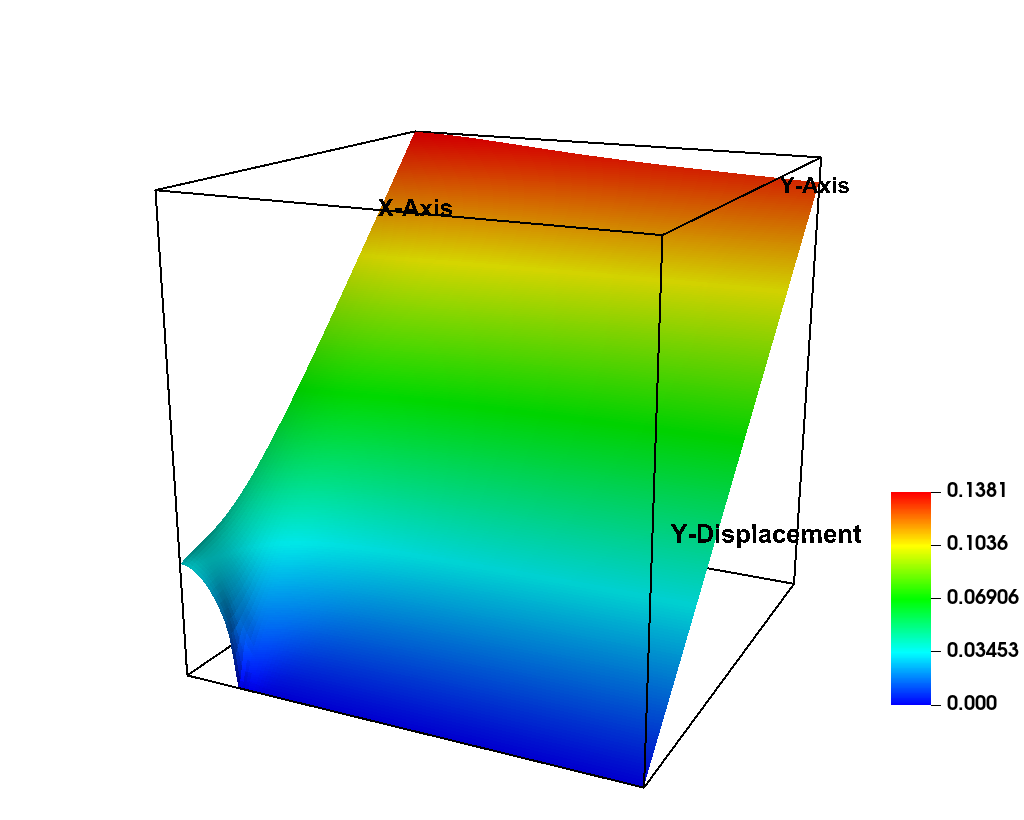}
        \caption{Y-displacement at $\alpha = 1.0$, $\sigma_{T} = 0.1$, and $\beta = 10.0$ (3D view)}
        \label{fig:y_displacement_3d}
    \end{subfigure}

    \caption{Displacement plots for $\alpha = 1.0$, $\sigma_{T} = 0.1$, and $\beta = 10.0$.}
    \label{fig:displacement_combined}
\end{figure}

Figure~\ref{fig:displacement_combined} illustrates both displacement metrics, including a vector plot of the displacement and a three-dimensional elevated representation of the $y$-displacement. It is evident that, even in scenarios where the fiber orientations are orthogonal to the crack front, the vertical displacement of the crack face displays an elliptical profile similar to that observed in previous cases. This observation suggests a smooth, curved deformation of the crack faces rather than a sharp, angular opening, thereby implying a more gradual stress distribution around the crack tip.

\section{Conclusion}\label{conclusion}

This study tackled the development of a continuous Galerkin finite element method for approximating solutions to vector-valued quasi-linear elliptic boundary value problems, explicitly arising from modeling geometrically linear, transversely isotropic elastic materials with algebraically nonlinear constitutive laws. The constitutive response was assumed to be monotonic and Lipschitz continuous, ensuring the well-posedness of our continuous Galerkin formulation. Furthermore, we established the existence and uniqueness of the discrete solution by applying the Riesz representation theory. Building upon this foundation, we successfully demonstrated the accurate computation of crack-tip fields within transversely isotropic elastic solids. Using our continuous Galerkin method with a suitable nonlinear constitutive model, we effectively captured the complex, highly localized stress and strain behaviors inherent to crack singularities. The numerical results validated the efficacy of our approach in addressing the challenges posed by material anisotropy and nonlinear responses, particularly in the critical near-tip region. The proven well-posedness of the continuous formulation, combined with the established existence and uniqueness of the discrete solution, reinforces the robustness of our developed methodology.

This research offers a valuable tool for analyzing fracture mechanics problems in transversely isotropic materials, contributing to a deeper understanding of their structural integrity and failure mechanisms. Our nonlinear model's observed strain energy density behavior closely resembled LEFM's. This intriguing similarity suggests the potential applicability of LEFM's local fracture criteria for studying crack evolution within the framework of the nonlinear strain-limiting constitutive relations presented herein. Future work could investigate adaptive mesh refinement techniques and extend the analysis to dynamic fracture scenarios, further enhancing this study's scope and applicability.  
 
\section*{Acknoledgement}
The first author, SG, would like to thank the University of Texas Rio Grande Valley for providing a Presidential Research Fellowship during his PhD studies. SMM's work is supported by the National Science Foundation under Grant No. 2316905.

\bibliographystyle{plain}
\bibliography{ref}

\end{document}